\documentclass[11pt]{article}
\usepackage{amssymb}
\usepackage{times}

\title{On minimal non-potentially closed subsets of the plane.\indent}
\author{Dominique LECOMTE}
\date{\it ~Topology Appl.\rm ~154,1 (2007), 241-262}

\newcommand{\Ana}{{\it\Sigma}^{1}_{1}}

\newcommand{\Borel}{{\it\Delta}^{1}_{1}}

\newcommand{\ana}{{\bf\Sigma}^{1}_{1}}

\newcommand{\borel}{{\bf\Delta}^{1}_{1}}
\newcommand{\boraone}{{\bf\Sigma}^{0}_{1}}

\newcommand{\boraxi}{{\bf\Sigma}^{0}_{\xi}}

\newcommand{\borone}{{\bf\Delta}^{0}_{1}}
\newcommand{\borthree}{{\bf\Delta}^{0}_{3}}

\newcommand{\bormone}{{\bf\Pi}^{0}_{1}}
\newcommand{\bormtwo}{{\bf\Pi}^{0}_{2}}

\newcommand{\bormxi}{{\bf\Pi}^{0}_{\xi}}

\newtheorem{thm} {Theorem}
\newtheorem{defi} [thm] {Definition}
\newtheorem{cor} [thm] {Corollary}
\newtheorem{lem} [thm] {Lemma}

\newtheorem{prop} [thm] {Proposition}

\begin{document}

\maketitle

\noindent {\footnotesize {\bf Abstract.} We study the Borel subsets of the plane that can be made closed by 
refining the Polish topology on the real line. These sets are called potentially 
closed. We first compare Borel subsets of the plane using products of 
continuous functions. We show the existence of a perfect antichain made of 
minimal sets among non-potentially closed sets. We apply this result to graphs, 
quasi-orders and partial orders. We also give a non-potentially closed set 
minimum for another notion of comparison. Finally, we show that we cannot have 
injectivity in the Kechris-Solecki-Todor\v cevi\' c dichotomy about analytic 
graphs.}\bigskip\bigskip

\noindent\bf {\Large 1 Introduction.}\rm\bigskip

 The reader should see [K] for the descriptive set theoretic notation used in 
this paper. This work is the continuation of a study made in [L1]-[L4]. 
The usual way of comparing Borel equivalence relations $E\subseteq 
X\times X$ and ${E'\subseteq X'\times X'}$ on Polish spaces is the 
Borel reducibility quasi-order:
$$E\leq_B E'~\Leftrightarrow ~\exists u\!:\! X\!\rightarrow\! 
X'~\hbox{\rm Borel~with}~E\! =\! (u\!\times\! u)^{-1}(E')$$
(recall that a quasi-order is a reflexive and transitive relation). 
Note that this makes sense even if $E$, $E'$ are not equivalence 
relations. It is known 
that if $(B_n)$ is a sequence of Borel subsets of $X$, then we can find a finer Polish 
topology on $X$ making the $B_n$'s clopen (see Exercise 13.5 in [K]). So assume 
that $E\leq_B E'$ and let $\sigma$ be a finer Polish topology on $X$ 
making $u$ continuous. If $E'$ is in some Baire class $\Gamma$, then 
$E\in\Gamma ([X,\sigma ]^2)$. This motivates the following (see [Lo2]):

\begin{defi} (Louveau) Let $X,Y$ be Polish spaces, 
$A$ a Borel subset of $X \times Y$, and $\Gamma$ a Baire class. We 
say that $A$ is $potentially~in~\Gamma$ 
$\big($denoted $A\in \mbox{pot}(\Gamma)\big)$ iff we can find a finer Polish 
topology $\sigma$ $($resp., $\tau )$ on $X$ $($resp., $Y)$ such that 
$A\in\Gamma ([X, \sigma ]\times [Y, \tau ])$.\end{defi}

 This notion is a natural invariant for $\leq_B$: if $E'\in \mbox{pot}(\Gamma)$ and 
$E\leq_B E'$, then $E\in \mbox{pot}(\Gamma)$. Using this notion, 
A. Louveau showed that the collection of $\boraxi$ equivalence 
relations is not cofinal for $\leq_B$, and deduces from this the 
inexistence of 
a maximum Borel equivalence relation for $\leq_B$.\bigskip 

 A. Louveau has also more recently noticed that one can associate a 
quasi-order $R_{A}\subseteq (X\times 2)^2$ to $A\subseteq X^2$ as follows: 
$$(x,i)~R_{A}~(y,j)~~\Leftrightarrow ~~(x,i)=(y,j)~~\mbox{or}~~
[(x,y)\in A~~\mbox{and}~~(i,j)=(0,1)].$$ 
Note that $R_{A}$ is also antisymmetric, so that it is actually a partial order. 

\vfill\eject

 A. Louveau noticed the following facts, using the following notion of 
comparison between Borel subsets $A\!\subseteq\! X\times Y$, 
${A'\!\subseteq\! X'\!\!\times\! Y'}$ of products of two Polish spaces:
$$A\sqsubseteq^r_{B} A'~\Leftrightarrow ~\exists u\!:\! X\!\rightarrow\! X'\ 
\exists v\!:\! Y\!\rightarrow\! Y'~\hbox{\rm one-to-one\ Borel~with}~A\! =\! 
(u\!\times\! v)^{-1}(A').$$
Here the letter $r$ means ``rectangle" ($u$ and $v$ may be different).\bigskip

\noindent - Assume that $A\subseteq X^2$ has full projections, and that 
$A'\subseteq (X')^2$. Then $A\sqsubseteq^r_{B}A'$ is equivalent to 
$R_{A}\leq_B R_{A'}$.\bigskip

\noindent - If $A\subseteq X^{2}$ is $\sqsubseteq^r_{B}$-minimal among 
non-potentially closed sets, then $R_{A}$ is $\leq_{B}$-minimal among 
non-potentially closed partial orders.\bigskip 

\noindent - Conversely, if $R_{A}$ is $\leq_{B}$-minimal among 
non-potentially closed partial orders and if $A$ has full projections, then $A$ 
is $\sqsubseteq^r_{B}$-minimal among non-potentially closed sets.\bigskip 

 These facts show that, from the point of view of Borel reducibility, the study 
of Borel partial orders is essentially the study of arbitrary Borel subsets 
of the plane. This strengthens the motivation for studying arbitrary Borel 
subsets of the plane, from the point of view of potential complexity.\bigskip

\noindent $\bullet$ A standard way to see that a set is complicated 
is to notice that it is more complicated than a well-known example. 
For instance, we have the following result (see [SR]):

\begin{thm} (Hurewicz) Let $P_{f}:=\{\alpha\!\in\! 
2^\omega /\exists n\!\in\!\omega ~\forall m\!\geq\! n~~\alpha (m)\! 
=\! 0\}$, $X$ be a Polish space and $A$ a Borel subset of $X$. Then 
exactly one of the following holds:\smallskip

\noindent (a) The set $A$ is $\bormtwo (X)$.\smallskip

\noindent (b) There is a continuous injection 
$u:2^\omega\rightarrow X$ such that $P_{f}=u^{-1}(A)$.
\end{thm}

 This theorem has been generalized to all Baire classes in [Lo-SR]. 
We try to adapt this result to the Borel subsets of the plane. In 
this direction, we have the following result for equivalence 
relations (see [H-K-Lo]):

\begin{thm} (Harrington-Kechris-Louveau) Let $X$ be a 
Polish space, $E$ a Borel equivalence relation on $X$, and 
$E_{0}\!:=\!\{(\alpha ,\beta )\!\in\! 2^\omega\times 2^\omega 
/\exists n\!\in\!\omega ~\forall m\!\geq\! n~~\alpha (m)\! =\!\beta 
(m)\}$. Then exactly one of the following holds:\smallskip

\noindent (a) The relation $E$ is $\mbox{pot}(\bormone)$.\smallskip

\noindent (b) We have $E_{0}\leq_B E$ (with $u$ continuous and 
one-to-one).\end{thm}

\noindent $\bullet$ We will study structures other than equivalence 
relations (for example quasi-orders), and even arbitrary Borel 
subsets of the plane. We need some other notions of comparison.\bigskip 

 Recall that Wadge's quasi-order $\leq_{W}$ on Borel subsets of $\omega^\omega$ 
is defined by
$$A\leq_W A'~~\Leftrightarrow ~~\exists 
u\!:\!\omega^\omega\!\rightarrow\! 
\omega^\omega ~\hbox{\rm continuous~with}~A\! =\! u^{-1}(A').$$
It is known that this quasi-order is well-founded (in the sense that there is no 
sequence $(B_{n})$ with $B_{n+1}\leq_{W}B_{n}$ and 
$B_{n}\not\leq_{W}B_{n+1}$ 
for each $n$). Moreover, any $\leq_{W}$-antichain is of cardinality 
at most 2 (in fact of the form $\{A,\neg A\}$). It follows that any class 
$\borel\setminus\bormxi$ admits a unique (up to the equivalence 
associated to $\leq_{W}$) minimal element.

\vfill\eject

\noindent $\bullet$ There are several natural ways of comparing Borel 
subsets 
$A\!\subseteq\! X\times Y$, ${A'\!\subseteq\! X'\!\!\times\! Y'}$ of 
products of two Polish spaces. All of them will have the same behavior here. 
The one we 
will use is the following:
$$A\leq^r_{c} A'~\Leftrightarrow ~\exists u\!:\! X\!\rightarrow\! 
X'\ \exists 
v\!:\! Y\!\rightarrow\! Y'~\hbox{\rm continuous~with}~A\! =\! 
(u\!\times\! v)^{-1}(A').$$
Here the letter $c$ is for ``continuous". We have the following (see [L1]):

\begin{thm} Let $\Delta (2^\omega )\!:=\!
\{(\alpha ,\beta )\!\in\! 2^\omega\!\times 2^\omega\! /\alpha\! 
=\!\beta\}$,  
${L_0\!:=\!\{(\alpha ,\beta )\!\in\! 2^\omega\!\!\times\! 2^\omega\!  
/\alpha\! <_{\mbox{lex}}\!\beta\}}$, $X$ and $Y$ be Polish spaces, and $A$ a 
$pot\big(\check D_2(\boraone)\big)$ subset of $X\times Y$. 
Then exactly one of the following holds:\smallskip

\noindent (a) The set $A$ is $\mbox{pot}(\bormone)$.\smallskip

\noindent (b) $\neg\Delta (2^\omega )\leq^r_{c} A$ or 
$L_{0}\leq^r_{c} A$ (with $u$ and $v$ one-to-one).\end{thm}

\noindent $\bullet$ Things become much more complicated at the level 
$D_{2}(\boraone)$ 
(differences of two open sets; $\check D_2(\boraone)$ is the dual 
Wadge class of unions of a closed set and an open set; notice that 
we can extend Definition 1 to the class $\check D_2(\boraone)$). We 
will show the following:

\begin{thm}  There is a perfect 
$\leq^r_{c}$-antichain 
$(A_{\alpha})_{\alpha\in 2^\omega}\subseteq D_{2}(\boraone)
(2^\omega\times 2^\omega )$ such that $A_{\alpha}$ is 
$\leq^r_{c}$-minimal among 
$\borel\setminus \mbox{pot}(\bormone)$ sets, for any $\alpha\!\in\! 2^\omega$.
\end{thm}

 In particular, unlike for classical Baire classes and $\leq_{W}$, 
one cannot characterize 
non-$\mbox{pot}(\bormone )$ sets by an obstruction condition involving only 
one (or even countably many) 
set(s). We will also show that ${[D_{2}(\boraone)\setminus 
\mbox{pot}(\bormone),\leq^r_{c}]}$ is ill-founded.\bigskip

 Theorem 5 can be applied to structures. We will show the following:

\begin{thm} There is a perfect 
$\leq_{B}$-antichain 
$(R_{\alpha})_{\alpha\in 2^\omega}\subseteq D_{2}(\boraone)
\big((2^\omega\times 2)^2\big)$ such that $R_{\alpha}$ 
is $\leq_{B}$-minimal among $\borel\setminus \mbox{pot}(\bormone)$ sets, for 
any 
$\alpha\in 2^\omega$. Moreover, $(R_{\alpha})_{\alpha\in 2^\omega}$ 
can be taken to be a subclass of any of the following classes:\smallskip

\noindent - Directed graphs (i.e., irreflexive relations).\smallskip

\noindent - Graphs (i.e., irreflexive and symmetric relations).\smallskip

\noindent - Oriented graphs (i.e., irreflexive and antisymmetric relations).\smallskip

\noindent - Quasi-orders.\smallskip

\noindent - Strict quasi-orders (i.e., irreflexive and transitive relations).\smallskip

\noindent - Partial orders.\smallskip

\noindent - Strict partial orders (i.e., irreflexive, antisymmetric and transitive 
relations).\end{thm}

\noindent\bf Remarks.\rm ~(a) Theorem 6 shows that Harrington, Kechris 
and Louveau's Theorem is very specific, and that the combination of 
symmetry and 
transitivity is very strong.\bigskip

\noindent (b) We produce concrete examples of such antichains. These examples 
must be in any complete family of minimal sets, up to bi-reducibility.

\vfill\eject

\noindent $\bullet$ Theorem 5 shows that any complete family of 
minimal sets for $[\borel\setminus \mbox{pot}(\bormone),\leq^r_{c}]$ has 
size continuum. So we must find another notion of comparison. In 
{[K-S-T]}, the following notion is defined. Let $X$, $X'$ be Polish 
spaces, and $A\!\subseteq\! X\!\times\! X$, $A'\!\subseteq 
X'\!\times\! X'$ be analytic sets. We set
$$(X,A)\preceq_{c}(X',A')\Leftrightarrow\exists u\!:\! 
X\!\rightarrow\! X'~\hbox{\rm continuous~with}~A\!\subseteq\! 
(u\!\times\! u)^{-1}(A').$$
 When $u$ is Borel we write $\preceq_{B}$ instead of $\preceq_{c}$.\bigskip
 
 Let $\psi:\omega\rightarrow 2^{<\omega}$ be the natural bijection 
($\psi (0)=\emptyset$, $\psi (1)=0$, $\psi (2)=1$, $\psi (3)=0^2$, 
$\psi (4)=01$, $\psi (5)=10$, $\psi (6)=1^2$, $\ldots$). Note that 
$|\psi (n)|\leq n$, so that we can define $s_{n}:=\psi (n)0^{n-|\psi 
(n)|}$. The 
crucial properties of $(s_{n})$ are that it is $dense$ (there is $n$ 
such that $t\prec s_{n}$, for each $t\in 2^{<\omega}$), and that 
$|s_{n}|=n$. We set
$$A_{1}\!:=\!\{ (s_{n}0\gamma ,s_{n}1\gamma )/n\!\in\!\omega 
~~\mbox{and}~\ \gamma\!\in\! 2^\omega\}.$$
The symmetric set $s(A_1)$ generated by $A_1$ is considered in 
[K-S-T], where the following is essentially proved:

\begin{thm} (Kechris, Solecki, Todor\v cevi\' c) 
Let $X$ be a Polish space and $A$ an analytic subset of $X\times X$. 
Then exactly one of the following holds:\smallskip

\noindent (a) $(X,A)\preceq_{B}(\omega ,\not=)$.\smallskip

\noindent (b) $(2^\omega ,A_{1})\preceq_{c}(X,A)$.\end{thm}

 Actually, the original statement in [K-S-T] is when $A$ is a graph, 
and with $s(A_1)$ instead of $A_1$. But we can get Theorem 7 without 
any change in the proof in [K-S-T].\bigskip 

\noindent $\bullet$ In [L3] the following is shown (see Theorem 2.9):

\begin{thm} Let $X$, $Y$ be Polish spaces, and $A$ a 
$\mbox{pot}(\borthree )$ subset of $X \times Y$. Then exactly one of the 
following holds:\smallskip

\noindent (a) The set $A$ is $\mbox{pot}(\bormone )$.\smallskip

\noindent (b) There are $u:2^\omega\rightarrow X$, $v:2^\omega\rightarrow Y$ 
continuous with ${A_{1}\! =\! (u\!\times\! v)^{-1}(A)\cap\overline{A_{1}}}$.\end{thm} 

\noindent (We can replace $A_1$ in [L3] by what we call $A_{1}$ 
here.) We generalize this result to arbitrary Borel subsets of $X\times Y$:

\begin{thm} Let $X$, $Y$ be Polish spaces, and $A$, 
$B$ be disjoint analytic subsets of 
$X \times Y$. Then exactly one of the following 
holds:\smallskip

\noindent (a) The set $A$ is separable from $B$ by a $\mbox{pot}(\bormone )$ set.\smallskip

\noindent (b) There are $u:2^\omega\rightarrow X$ and 
$v:2^\omega\rightarrow Y$ 
continuous such that the inclusions $A_{1}\subseteq  (u\times v)^{-1}(A)$ and 
$\overline{A_{1}}\setminus A_1\subseteq (u\times v)^{-1}(B)$ hold.\smallskip

\noindent Moreover, we can neither replace $\overline{A_{1}}\setminus A_1$ 
with $(2^\omega\times 2^\omega)\setminus A_1$, nor ensure that 
$u$ and $v$ are one-to-one.\end{thm} 

 So we get a minimum non-potentially closed set if we do not ask for 
a reduction on the whole product.

\vfill\eject

\noindent $\bullet$ In [K-S-T], it is conjectured that we can have 
$u$ one-to-one in Theorem 7.(b). This is not the case:

\begin{thm} There is no graph $(X_{0},R_{0})$ with 
$X_{0}$ Polish and $R_{0}\in\ana (X_{0}^2)$ such that for every graph 
$(X,A)$ of the same type, exactly one of the following 
holds:\smallskip

\noindent (a)~$(X,A)\preceq_{B}(\omega ,\not= )$.\smallskip

\noindent (b)~$(X_{0},R_{0})\preceq_{c,1-1}(X,A)$.\end{thm}

  The proof is based on the counterexample used in [L3] to show that 
we cannot have injectivity in Theorem 2.9.\bigskip 

\noindent $\bullet$ The paper is organized as follows.\bigskip

\noindent - In Section 2, we prove Theorem 9.\bigskip

\noindent - In Section 3, we prove Theorem 10.\bigskip

\noindent - In Section 4, we give a sufficient condition for 
minimality among non-potentially closed sets. We use it to prove 
Theorems 5 and 6.\bigskip

\noindent - In Section 5, we give conditions on $A$ which allow us to 
replace $\overline{A_{1}}\setminus A_1$ with $(2^\omega\times 
2^\omega)\setminus A_1$ in Theorem 9 (and therefore come back to 
$\leq^r_{c}$). We can write 
${A_{1}\! =\!\bigcup_{n}~\mbox{Gr}(f_{n})}$, where ${f_{n}(s_{n}0\gamma 
)\!:=\! s_{n}1\gamma}$. 
Roughly speaking, we require that the $f_{n}$'s do not induce cycles. 
This is really the key property making the $A_{\alpha}$'s appearing 
in the statement of Theorem 5 pairwise orthogonal. We will deduce 
from this the minimality of $A_1$ among non-potentially closed sets 
for $\leq^r_{c}$, using the sufficient condition for minimality in Section 4.\bigskip\smallskip

\noindent\bf {\Large 2 A minimum non-potentially closed set.}\rm\bigskip

 We will prove Theorem 9. The proof illustrates the link between the 
dichotomy results in [K-S-T] and the notion of potential Baire class. 
We will see another link in Section 3. The next lemma is essentially 
Lemma 3.5 in [L1], and the crucial point of its proof.

\begin{lem} Let $X$ be a nonempty Polish space, $n$ be 
an integer, $D_{f_n}$ and $f_n[D_{f_n}]$ be dense $G_{\delta}$ 
subsets of some open subsets of $X$, and 
$f_{n}:D_{f_n}\rightarrow f_n[D_{f_n}]$ a continuous and open 
map.\smallskip

\noindent (a) Let $G$ be a dense $G_{\delta}$ subset of $X$. Then 
$\mbox{Gr}(f_{n})\subseteq\overline{\mbox{Gr}(f_{n})\cap G^{2}}$, for each $n$.\smallskip

\noindent (b) Let $A:=\bigcup_{n}~\mbox{Gr}(f_{n})$. If $\Delta 
(X)\subseteq\overline{A}\setminus A$, then $A$ is not $\mbox{pot}(\bormone )$.\end{lem}

\noindent {\bf Proof.} (a) Let $U$ (resp., $V$) be an open 
neighborhood of 
$x\in D_{f_{n}}$ (resp., $f_{n}(x)$). Then $f_{n}[D_{f_{n}}]\cap 
V\cap G$ is a 
dense $G_{\delta}$ subset of $f_{n}[D_{f_{n}}]\cap V$, thus 
$f_n^{-1}(V\cap G)$ 
is a dense $G_{\delta}$ subset of $f_n^{-1}(V)$. Therefore $G\cap 
f_n^{-1}(V)$ 
and $G\cap f_n^{-1}(V\cap G)$ are dense $G_{\delta}$ subsets of 
$f_n^{-1}(V)$. So we can find 
$$y\in U\cap G\cap f_n^{-1}(V\cap G).$$ 
Now $\big(y,f_n(y)\big)$ is in 
the intersection $(U\times V)\cap \hbox{Gr}(f_n)\cap 
G^2$, so this set is non-empty.\bigskip

\noindent (b) We argue by contradiction: we can find a finer Polish topology on 
$X$ such that $A$ becomes closed. By 15.2, 11.5 and 8.38 in [K], the new topology 
and the old one agree on a dense $G_{\delta}$ subset of $X$, say $G$: 
$A\cap G^2\in\bormone (G^2)$. Let $x\in G$. We have 
$(x,x)\in G^2\cap\overline{A}\setminus A$. By (a) we get 
$\overline{A}\subseteq\overline{A\cap G^2}$. Thus 
$(x,x)\in G^2\cap\overline{A\cap G^2}\setminus (A\cap G^{2})$, which 
is absurd.\hfill{$\square$} 

\vfill\eject

\begin{cor}  The set 
$A_{1}=\overline{A_{1}}\setminus\Delta (2^\omega )$ is 
$D_{2}(\boraone )\setminus \mbox{pot}(\bormone )$, and 
$\overline{A_{1}}=A_{1}\cup\Delta (2^\omega )$.\end{cor}

\noindent {\bf Proof.} As we saw in the introduction, we can write 
${A_{1}\! =\!\bigcup_{n}~\hbox{Gr}(f_{n})}$, where 
${f_{n}(s_{n}0\gamma )\!:=\! s_{n}1\gamma}$. Notice that $f_{n}$ is a 
partial homeomorphism with clopen domain and range. Moreover, we have 
$$\Delta (2^\omega )\subseteq\overline{A_{1}}\setminus A_{1}$$ 
(in fact, the equality holds). 
Indeed, if $t\in 2^{<\omega}$, we have $(s_{\psi^{-1}(t)}0^\infty ,
s_{\psi^{-1}(t)}10^\infty )\in N^2_{t}\cap A_{1}$. Thus 
$A_{1}=\overline{A_{1}}\setminus\Delta (2^\omega )$ is 
$D_{2}(\boraone )$, and the corollary follows from Lemma 
11.\hfill{$\square$}\bigskip 

\noindent {\bf Proof of Theorem 9.} We cannot have (a) and (b) 
simultaneously. For if $D$ is 
potentially closed and separates $A$ from $B$, then we get 
$A_1=(u\times v)^{-1}(D)\cap\overline{A_1}$, thus $A_{1}\in 
\mbox{pot}(\bormone )$, which contradicts Corollary 12.\bigskip  

\noindent $\bullet$ Let $f:\omega^\omega\rightarrow X\times Y$ be 
a continuous map with 
$f[\omega^\omega ]=B$, and $f_{0}$ (resp., $f_{1}$) be the first 
(resp., second) coordinate of $f$, so that 
${(f_{0}\times f_{1})[\Delta (\omega^\omega )]=B}$. We set 
$R:=(f_{0}\times f_{1})^{-1}(A)$, which is an irreflexive analytic 
relation on $\omega^\omega$. By Theorem 7, either there exists a 
Borel map ${c:\omega^\omega\rightarrow\omega}$ such that 
$(\alpha ,\beta )\in R$ implies ${c(\alpha )\not= c(\beta )}$, or there 
is a continuous map $u_{0}:2^\omega\rightarrow\omega^\omega$ such that 
${(\alpha ,\beta )\in A_{1}}$ implies 
${\big(u_{0}(\alpha ),u_{0}(\beta )\big)\in R}$.\bigskip 

\noindent $\bullet$ In the first case, we define ${C_{n}:=c^{-1}(\{ n\})}$. We get 
${\Delta (\omega^\omega )\subseteq\bigcup_{n} C_{n}^2\subseteq\neg R}$, so that 
$${B\subseteq\bigcup_{n} f_{0}[C_{n}]\times f_{1}[C_{n}]\subseteq\neg A}.$$ 
By a standard reflection argument there is a sequence $(X_n)$ (resp., $(Y_n)$) of Borel subsets of $X$ (resp., $Y$) with  
$$\bigcup_{n} f_{0}[C_{n}]\times f_{1}[C_{n}]\subseteq\bigcup_{n} X_{n}\times Y_n\subseteq\neg A.$$
But $\bigcup_{n}\ X_{n}\times Y_n$ is $\mbox{pot}(\boraone)$, so we are in the case (a).\bigskip

\noindent $\bullet$ In the second case, let $u:=f_{0}\circ u_{0}$, 
$v:=f_{1}\circ u_{0}$. These maps satisfy the conclusion of condition (b) because 
$\overline{A_{1}}\setminus A_{1}\subseteq\Delta (2^\omega )$, by Corollary 12.\bigskip

\noindent $\bullet$ By the results in [L3], we can neither replace 
$\overline{A_{1}}\setminus A_1$ with $(2^\omega\times 
2^\omega)\setminus A_1$, nor can we ensure that $u$ and $v$ are 
one-to-one.\hfill{$\square$}\bigskip

\noindent\bf Remarks.\rm ~(a) In Theorem 9, we cannot ensure that $u=v$ when 
$X=Y$: take $X:=2^\omega$, 
$$A:=\{ (\alpha ,\beta)\in N_{0}\times N_{1}/\alpha <_{\mbox{lex}}\beta\}$$ 
and $B:=(N_{0}\times N_{1})\setminus A$.\bigskip

\noindent (b) This proof cannot be generalized, in the sense that we used the 
fact that the range of a countable union of Borel rectangles (a $\mbox{pot}(\boraone)$ 
set) by a product function is still a countable union of rectangles, so more or 
less a $\mbox{pot}(\boraone)$ set. This fails completely for the dual level. Indeed, we 
saw that the range of the diagonal (which is closed) by a product function can be 
any analytic set. So in view of generalizations, it is better to have another 
proof of Theorem 9.

\vfill\eject

\noindent\bf {\Large 3 The non-injectivity in the Kechris-Solecki-Todor\v cevi\' c 
dichotomy.}\rm\bigskip 

 Now we will prove Theorem 10. The proof we give is not the original one, which 
used effective descriptive set theory, and a reflection argument. The proof we 
give here is due to B. D. Miller, and is a simplification of the original proof.
\bigskip

\noindent\bf Notation.\rm ~If $A\subseteq X^2$, 
$A^{-1}:=\{(y,x)\in X^2/(x,y)\in A\}$ and $s(A):=A\cup A^{-1}$ is the symmetric 
set generated by $A$.\bigskip

\noindent $\bullet$ Fix sets $S_{0}\supseteq S_{1}\supseteq\ldots$ of 
natural numbers such that\bigskip

(1) $S_{n}\setminus S_{n+1}$ is infinite for each integer $n$.\bigskip

(2) $\bigcap_{n\in\omega}\ S_{n}=\emptyset$.\bigskip

\noindent $\bullet$ For each $n\in\omega$, fix 
$f_{n}:S_{n}\rightarrow S_{n}\setminus S_{n+1}$ injective, and define 
$g_{n}:2^\omega\rightarrow 2^\omega$ by 
$$[g_{n}(\alpha )](k):=\left\{\!\!\!\!\!\!\!\!
\begin{array}{ll} 
& \alpha [f_{n}(k)]~\hbox{\rm if}~k\in S_{n}\mbox{,}\cr\cr 
& \alpha (k)~\hbox{\rm otherwise.}
\end{array}
\right.$$
$\bullet$ It is clear that each of the closed sets 
$M_{n}:=\{ \alpha \in 2^\omega /g_{n}(\alpha )=\alpha\}$ is meager, and since 
each $g_{n}$ is continuous and open, it follows that the $F_{\sigma}$ set 
$$M:=\bigcup_{s\in\omega^{<\omega}, n\in\omega}\ (g_{s(0)}\circ\ldots\circ 
g_{s(\vert s\vert -1)})^{-1}(M_{n})$$
is also meager, so that $X:=2^\omega\setminus M$ is a comeager, dense 
$G_{\delta}$ set which is invariant with respect to each $g_{n}$. 
Put $G_{1}:=\bigcup_{n\in\omega}\ s[\hbox{Gr}({g_{n}}\vert_{X})]$.\bigskip

\noindent {\bf Proof of Theorem 10.} We argue by contradiction: this gives 
$(X_{0},R_{0})$.\bigskip

\noindent {\bf Claim\ 1.}\it\ Let $X$ be a Polish space, 
and $g_{0},g_{1},\ldots :X\rightarrow X$ fixed-point free Borel functions 
such that $g_{m}\circ g_{n}=g_{m}$ if $m<n$. Then every locally countable Borel 
directed subgraph of the Borel directed graph 
$G:=\bigcup_{n\in\omega}\ \hbox{Gr}(g_{n})$ has countable Borel chromatic number, 
i.e., satisfies Condition (a) in Theorem 7.\bigskip\rm 

 Suppose that $H$ is a locally countable Borel directed subgraph of $G$. By the 
Lusin-Novikov uniformization theorem, there are Borel partial injections $h_{n}$ 
on $X$ such that $H=\bigcup_{n\in\omega}\ \hbox{Gr}(h_{n})$. By replacing each 
$h_{n}$ with its restrictions to the sets $\{ x\in D_{h_{n}}/h_{n}(x)=g_{m}(x)\}$, 
for $m\in\omega$, we can assume that for all $n\in\omega$, there is 
$k_{n}\in\omega$ such that $h_{n}=g_{k_{n}}\vert_{D_{h_{n}}}$. It is easily seen 
that the directed graph associated with a Borel function has countable Borel 
chromatic number (see also Proposition 4.5 of [K-S-T]), so by replacing $h_{n}$ 
with its restriction to countably many Borel sets, we can assume also that for 
all $n\in\omega$, 
$D_{h_{n}}^{2}\cap\bigcup_{k\leq k_{n}}\ \hbox{Gr}(g_{k})=\emptyset$. It only 
remains to note that 
$D_{h_{n}}^{2}\cap\bigcup_{k>k_{n}}\ \hbox{Gr}(g_{k})=\emptyset$. To see this, 
simply observe that if $k>k_{n}$ and $x,g_{k}(x)\in D_{h_{n}}$, then 
$h_{n}(x)=g_{k_{n}}(x)=g_{k_{n}}\circ g_{k}(x)=h_{n}\circ g_{k}(x)$, which 
contradicts the fact that $h_{n}$ is a partial injection. This proves the claim.\hfill{$\diamond$}

\vfill\eject

\noindent {\bf Claim\ 2.}\it\ The Borel graph $G_{1}$ has uncountable Borel chromatic 
number, but if $H\subseteq G_{1}$ is a locally countable Borel directed graph, 
then $H$ has countable Borel chromatic number.\bigskip\rm

 Condition (1) implies that $g_{m}\circ g_{n}=g_{m}$ if $m<n$, so Claim 1 ensures 
that if $H\subseteq G_{1}$ is a locally countable Borel directed graph, then $H$ 
has countable Borel chromatic number.\bigskip

 To see that $G_{1}$ has uncountable Borel chromatic number, it is enough to show 
that if $B\in\borel (2^\omega )$ is non-meager, then 
$B\cap G_{1}^{2}\not=\emptyset$. Let $s\in 2^{<\omega}$ such that $B$ is comeager 
in $N_{s}$. It follows from condition (2) that there is $n\in\omega$ such that 
$\vert s\vert <k$ for each $k\in S_{n}$. Then $g_{n}$ is a continuous, open map 
which sends $N_{s}$ into itself, thus 
$B\cap X\cap N_{s}\cap g_{n}^{-1}(B\cap X\cap N_{s})$ is comeager in $N_{s}$. 
Letting $x$ be any element of this set, it follows that $x$, $g_{n}(x)$ are 
$G_{1}$-related elements of $B$.\hfill{$\diamond$}\bigskip

 We are now ready to prove the theorem: as $(X_{0},R_{0})$ satisfies (b), it does 
not satisfy (a). Therefore $R_{0}$ has uncountable Borel chromatic number. As 
$s(A_{1})$ and $G_{1}$ have uncountable Borel chromatic number, we get 
$(X_{0},R_{0})\preceq_{c,1-1}[2^\omega ,s(A_{1})]$ and 
$(X_{0},R_{0})\preceq_{c,1-1}(2^\omega ,G_{1})$ (with witness $\pi$). As 
$s(A_{1})$ is locally countable, $R_{0}$ is also locally countable. Therefore 
$(\pi\times\pi )[R_{0}]$ is a locally countable Borel subgraph of $G_{1}$ with 
uncountable Borel chromatic number, which contradicts Claim 2.\hfill{$\square$}
\bigskip

\noindent\bf Remark.\rm ~This proof also shows a similar theorem for irreflexive 
analytic relations, by considering 
$\bigcup_{n\in\omega}\ \hbox{\rm Gr}(g_{n}\vert_{X})$ (resp., $A_{1}$) 
instead of $G_{1}$ (resp., $s(A_{1})$).\smallskip\bigskip 

\noindent\bf {\Large 4 Perfect antichains made of sets minimal among 
non-$\mbox{pot}(\bormone )$ sets.}\bigskip\rm

 As mentioned in the introduction, a great variety of very different 
examples appear at level $D_{2}(\boraone)$, all of the same type. Let us make 
this more specific.

\begin{defi}  We say that $\big(X,(f_{n})\big)$ is a $converging~situation$ if\smallskip

\noindent (a) $X$ is a nonempty 0-dimensional perfect Polish space.\smallskip
  
\noindent (b) $f_{n}$ is a partial homeomorphism with $\borone(X)$ domain and 
range.\smallskip

\noindent (c) The diagonal $\Delta (X)=\overline{A^f}\setminus A^f$, where 
$A^f:=\bigcup_{n}~\mbox{Gr}(f_{n})$.\end{defi}

 This kind of situation plays an important role in the theory of 
potential complexity (see, for example, Definition 2.4 in [L3]).\bigskip

\noindent {\bf Remarks.}~(a) Note that if $\big(X,(f_{n})\big)$ is a 
converging situation, then Lemma 11 ensures that $A^f$ is  
$D_{2}(\boraone)\setminus \mbox{pot}(\bormone)$, since 
${A^f=\overline{A^f}\setminus\Delta (X)}$.\bigskip

 \noindent (b) It is clear that an analytic graph $(X,A)$ has  
countable Borel chromatic number if and only if $A$ is separable from 
$\Delta (X)$ by a $\mbox{pot}(\borone)$ set. By Remark (a), this implies that 
$\big(2^\omega ,s(A^f)\big)$ does not have countable Borel chromatic number if 
$\big(X,(f_{n})\big)$ is a converging situation.\bigskip

\noindent\bf Notation.\rm ~In the sequel, we set 
$f^{B}_{n}:={f_{n}\vert}_{B\cap f_{n}^{-1}(B)}$ if $B\subseteq X$ 
and $\big(X,(f_{n})\big)$ is a converging situation, so that 
$\hbox{Gr}(f^{B}_{n})=\hbox{Gr}(f_{n})\cap B^2$.\bigskip

 The reader should see [Mo] for the basic notions of effective descriptive set theory. 
Let $Z$ be a recursively presented Polish space.\bigskip 

\noindent $\bullet$ The topology ${\it\Delta}_{Z}$ is the topology on $Z$ 
generated by $\Borel (Z)$. This topology is Polish (see the proof of 
Theorem 3.4 in [Lo2]).\bigskip 

\noindent $\bullet$ The Gandy-Harrington topology ${\it\Sigma}_{Z}$ on $Z$ is 
generated by $\Ana (Z)$. Recall that 
$$\Omega_{Z}:=\{z\in Z/\omega_1^z =\omega_1^{\mbox{CK}}\}$$ 
is Borel and $\Ana$, and $[\Omega_{Z},{\it\Sigma}_{Z}]$ is a $0$-dimensional Polish space (in fact, the 
intersection of $\Omega_{Z}$ with any nonempty $\Ana$ set is a nonempty clopen 
subset of $[\Omega_{Z},{\it\Sigma}_{Z}]$-see [L1]).

\begin{lem} Let $\big(X,(f_{n})\big)$ be a converging 
situation, $P$ a Borel subset of $X$ such that $A^f\cap P^2$ is not 
$\mbox{pot}(\bormone )$, and $\sigma$ a finer Polish topology on $P$. 
Then we can find a Borel subset $S$ of $P$ and a topology 
$\tau$ on $S$ finer than $\sigma$ such that 
$\big([S,\tau],(f^{S}_{n})_{n}\big)$ 
is a converging situation.\end{lem}

\noindent {\bf Proof.} We may assume that $[P,\sigma ]$ is recursively presented 
and $f^{P}_{n}$, $A^f\cap P^2$ are $\Borel$. We set 
${D\!:=\!\{x\!\in\! P/x\!\in\! {\it\Delta }^1_1\}}$, and 
${S\!:=\!\{ x\!\in\! P/(x,x)\!\in\!\overline{A^{f}\cap P^2}^
{{\it\Delta}_{P}^2}\}\cap\Omega_{P}\!\setminus\! D}$. As 
${S\!\in\!\Ana}$, 
$[S,{\it\Sigma }_{P}]$ is a 0-dimensional perfect Polish space. We 
set ${E:=A^{f}\cap (P\setminus D)^{2}}$. Note that $D$ is countable. 
By Remark 2.1 in [L1], $E$ is not potentially closed since
$$A^{f}\cap P^{2}=[A^{f}\cap \big( (P\cap D)\times P\big)]\cup 
[A^{f}\cap\big( P\times (P\cap D)\big)]\cup E.$$ 
Therefore ${\overline{E}^{{\it\Delta}_{P}^2}\setminus E}$ is a nonempty 
subset of ${(P\setminus D)^{2}\cap\overline{A^{f}}\setminus 
A^{f}\subseteq\Delta (X)}$. Thus ${S\not=\emptyset}$. Note also that  
${(x,x)\!\in\!\overline{A^{f}\cap P^2}^{{\it\Delta}_{P}^2}\cap S^2\! 
=\!\overline{A^{f}\cap P^2}^{{\it\Sigma}_{P}^2}\cap S^2\! =\!
\overline{A^{f}\cap S^2}^{[S,{\it\Sigma}_{P} ]^2}}$ if $x\!\in\! S$. 
Conversely, we have ${\overline{A^{f}\cap S^2}^{[S,{\it\Sigma}_{P}]^2}
\!\setminus\!\big(A^{f}\cap S^2\big)\!\subseteq\! 
S^2\cap\overline{A^{f}}\!
\setminus\! A^{f}\!\subseteq\!\Delta (S)}$. We have proved that $S$ is 
a Borel subset of $P$ such that 
$\big([S,{\it\Sigma}_{P}],(f^{S}_{n})_{n}\big)$ is a 
converging situation.\hfill{$\square$}

\begin{thm} Let $Y$, $Y'$ be Polish spaces, 
$A\in\borel (Y\times Y')$, $\big(X,(f_{n})\big)$ a converging 
situation. We 
assume that $A\leq^r_{c} A^f$. Then exactly one of the following 
holds:\smallskip

\noindent (a) The set $A$ is $\mbox{pot}(\bormone)$.\smallskip

\noindent (b) We can find a Borel subset $B$ of $X$ and a finer topology $\tau$ on $B$ 
such that $\big([B,\tau],(f^{B}_{n})_{n}\big)$ is a converging situation and $A^f\cap B^{2}\leq^r_{c}A$.
\end{thm}

\noindent {\bf Proof.} Let $u$ and $v$ be continuous functions such 
that $A=(u\times v)^{-1}(A^f)$. We assume that $A$ is not potentially 
closed. By Theorem 9 we can find continuous maps $u':2^\omega\rightarrow Y$ 
and $v':2^\omega\rightarrow Y'$ such that 
$A_{1}=(u'\times v')^{-1}(A)\cap\overline{A_{1}}$. We set 
$H:=u\big[u'[2^\omega ]\big]$, $K:=v\big[v'[2^\omega ]\big]$ and 
$P:=H\cap K$. Then $H$, $K$ and $P$ are compact and $A^{f}\cap (H\times K)$ is 
not $\mbox{pot}(\bormone)$ since 
$$A_{1}=[(u~\circ ~u')\times (v~\circ ~v')]^{-1}\big(A^{f}\cap (H\times K)\big)\cap\overline{A_{1}}$$ 
(we have $A_{1}\notin \mbox{pot}(\bormone )$ by Corollary 12). Therefore $A^{f}\cap P^{2}$ is not $\mbox{pot}(\bormone)$, since
$$\begin{array}{ll} 
A^{f}\cap (H\times K)\!\!\!\! 
& = [A^{f}\cap\big( (H\setminus K)\times K\big)]\cup 
[A^{f}\cap\big( H\times (K\setminus H)\big)]\cup [A^{f}\cap P^{2}]\cr 
& = [\overline{A^{f}}\cap\big( (H\setminus K)\times K\big)]\cup 
[\overline{A^{f}}\cap\big( H\times (K\setminus H)\big)]\cup 
[A^{f}\cap P^{2}].
\end{array}$$
By Lemma 14 we can find a Borel subset $S$ of $P$ and a finer 
topology $\sigma$ on $S$ such that $\big([S,\sigma ],(f^{S}_{n})_{n}\big)$ is a 
converging situation.

\vfill\eject

 By the Jankov-von Neumann Theorem there is $f':S\rightarrow 
u^{-1}(S)$ (respectively, ${g':S\!\rightarrow\! v^{-1}(S)}$) 
Baire measurable such 
that $u\big(f'(x)\big)\! =\! x$ (respectively, ${v\big(g'(x)\big)\! 
=\! x}$), for 
each $x\in S$. Notice that $f'$ and $g'$ are one-to-one. Let $G$ be a 
dense $G_{\delta}$ subset of $S$ such that $f'\vert_{G}$ and $g'\vert_{G}$ are 
continuous. These functions are witnesses to the inequality 
$A^f\cap G^2\leq^r_{c}A$. By Lemma 11, we get 
${\hbox{Gr}(f^{S}_{n})\subseteq\overline{\hbox{Gr}(f^{S}_{n})\cap G^2}}$. 
Therefore ${\overline{A^f\cap S^2}=\overline{A^f\cap G^2}}$, $\Delta 
(G)=G^{2}\cap\overline{A^f
\cap G^2}\setminus (A^f\cap G^2)$, and $A^f\cap G^2$ is not 
$\mbox{pot}(\bormone)$ by Lemma 11.\bigskip 

 By Lemma 14 we can find a Borel subset $B$ of $G$, equipped with 
some topology $\tau$ finer than $\sigma$, such that $\big([B,\tau 
],(f^{B}_{n})_{n}\big)$ is a converging situation.\hfill{$\square$}

\begin{cor}  Let $\big(X,(f_{n})\big)$ be a 
converging situation. The following statements are equivalent:\smallskip

\noindent (a) $A^f$ is $\leq^r_{c}$-minimal among $\borel\setminus 
\mbox{pot}(\bormone)$ sets.\smallskip
  
\noindent (b) For any Borel subset $B$ of $X$ and any finer Polish topology 
$\tau$ on $B$, $A^f \leq^r_{c} A^f\cap B^2$ if $A^f\cap B^2$ is not 
$\mbox{pot}(\bormone)$.\smallskip

\noindent (c) For any Borel subset $B$ of $X$ and for each finer topology $\tau$ 
on $B$, $A^f \leq^r_{c} A^f\cap B^2$ if $\big([B,\tau],(f^{B}_{n})_{n}\big)$ is a 
converging situation.\end{cor} 

\noindent {\bf Proof.} (a) $\Rightarrow$ (b) and (b) $\Rightarrow$ (c) 
are obvious. So let us show that (c) $\Rightarrow$ (a). Let $Y$, $Y'$ be 
Polish spaces, $A\in\borel (Y\times Y')\setminus \mbox{pot}(\bormone)$. We assume that 
$A\leq^r_{c} A^f$. By Theorem 15 we get a Borel subset $B$ of $X$ and 
a finer topology $\tau$ on $B$ such that $\big([B,\tau],(f^{B}_{n})_{n}\big)$ 
is a converging situation and $A^f\cap B^{2}\leq^r_{c}A$. By (c) we get 
$A^f \leq^r_{c} A^f\cap B^2$. Therefore $A^f \leq^r_{c} 
A$.\hfill{$\square$}\bigskip  

 This is the sufficient condition for minimality that we mentioned 
in the introduction. The following definitions, notation and facts 
will be used here and in Section 5 to build the reduction functions 
in the minimality results that we want to show.

\begin{defi} Let $R$ be a relation on a set 
$E$.\smallskip
 
\noindent $\bullet$ An $R\! -\! path$ is a finite 
sequence $(e_{i})_{i\leq n}\!\subseteq\! E$ such that 
$(e_{i},e_{i+1})\!\in\! R$ for $i\! <\! n$.\smallskip
 
\noindent $\bullet$ We say that $E$ is $R\! -\! connected$ if there is an $R$-path $(e_{i})_{i\leq n}$ with $e_0=e$ and $e_n=e'$ for each $e,\ e'\in E$.\smallskip
 
\noindent $\bullet$ An $R\! -\! cycle$ is an $R$-path 
$(e_{i})_{i\leq n}$ such that $n\!\geq\! 3$ and
$$[0\!\leq\! i\!\not=\! j\!\leq\! n~\hbox{\rm and}~e_{i}\! =\! 
e_{j}]~\Leftrightarrow ~\{ i,j\}\! =\!\{ 0,n\}.$$
$\bullet$ We say that $R$ is $acyclic$ if there is no $R$-cycle.\end{defi}

 Recall that if $R$ is symmetric and acyclic, $e$, $e'\!\in\! E$ and 
$(e_{i})_{i\leq n}$ is an $R$-path with $e_{0}=e$ and $e_{n}=e'$, 
then we can find a unique $R$-path $p_{e,e'}:=(f_{j})_{j\leq m}$ without 
repetition with $f_{0}=e$ and $f_{m}=e'$. We will write $|p_{e,e'}|=m+1$.
\bigskip

\noindent {\bf Notation.}~Let $\Theta:= (\theta_n)\subseteq 
2^{<\omega}$ with $|\theta_n|=n$. We will use two examples of such 
$\Theta$'s: $\theta_n = 0^n$ and $\theta_n = s_n$ (where $s_n$ has 
been defined in the introduction to build $A_1$). We define a tree 
${\mathfrak  R}_\Theta$ on $2\times 2$:
$${\mathfrak  R}_\Theta:=\{ (e,e')\!\in\! (2\!\times\! 2)^{<\omega}/e\! 
=\! e'\ \ \hbox{\rm or}\ \ \exists n\!\in\!\omega\ \exists w\!\in\! 
2^{<\omega}\ \ (e,e')\! =\! (\theta_n0w,\theta_n1w)\}.$$
Recall that $s({\mathfrak  R}_\Theta )$ is the symmetric set generated by 
${\mathfrak  R}_\Theta $.

\begin{prop} (a) $\big( 2^n,s({\mathfrak  R}_\Theta )\big)$ is connected, for each $n\in\omega$.\smallskip

\noindent (b) The relation $s({\mathfrak  R}_\Theta )$ is acyclic.\smallskip

\noindent (c) If $e,\ e'\in 2^n$ and $l<n$ is maximal with $e(l)\not= 
e'(l)$, the coordinate $l$ is 
changed only once in $p_{e,e'}$, and the other changed coordinates 
are at a level less than $l$.\end{prop}

\noindent {\bf Proof.} (a) We argue by induction on $n$. As 
$(\emptyset )$ is an 
$s({\mathfrak  R}_\Theta )$-path from $\emptyset$ to $\emptyset$, the 
statement is true for $n=0$. Assume that it is true at the level $n$, and 
let $e$, $e'\in 2^{n+1}$. We can write $e=s\epsilon$ and 
$e'=s'\epsilon'$, where $s,t\in 2^n$ and $\epsilon ,\epsilon'\in 2$. 
If $\epsilon =\epsilon'$, then let $(f_i)_{i\leq m}$ be an 
$s({\mathfrak  R}_\Theta )$-path with $f_0=s$ and $f_m=s'$. Let $e_i:=f_i\epsilon$. 
Then $(e_i)_{i\leq m}$ is an $s({\mathfrak  R}_\Theta 
)$-path with $e_0=e$ and 
$e_m=e'$. If $\epsilon\not=\epsilon'$, then let $(f_i)_{i\leq m}$ be an 
$s({\mathfrak  R}_\Theta )$-path with $f_0=s$ and $f_m=\theta_n$, and 
$(g_j)_{j\leq p}$ be an 
$s({\mathfrak  R}_\Theta )$-path with $g_0=\theta_n$ and $g_p=s'$. We set 
$e_i:=f_i\epsilon$ if $i\leq m$, 
$g_{i-m-1}\epsilon'$ if $m<i\leq m+p+1$. Then $(e_i)_{i\leq m+p+1}$ 
is an $s({\mathfrak  R}_\Theta )$-path with $e_0=e$ and 
$e_{m+p+1}=e'$.\bigskip

\noindent (b) We argue by contradiction. Let $(e_i)_{i\leq n}$ be an 
$s({\mathfrak  R}_\Theta )$-cycle, $p>0$ be the common length of the 
$e_i$'s, and $l<p$ maximal such that the sequence $(e_i(l))_{i\leq 
n}$ is not constant. We can find $i_1$ minimal with $e_{i_1}(l)\not= 
e_{i_1+1}(l)$. We have 
${e_{i_1}(l)\! =\! e_0(l)\! =\! e_n(l)}$. We can find $i_2>i_1+1$ 
minimal with 
${e_{i_1+1}(l)\not= e_{i_2}(l)}$. Then ${e_{i_1}(l)=e_{i_2}(l)}$ and 
$e_{i_1}=e_{i_2}$, because $|\theta_l|=l$. Thus $i_1=0$ and 
$i_2=n$. But $e_{i_1+1} = e_{i_2-1}$, which is absurd. Note that this 
proof of (b) is essentially in [L3], Theorem 2.7.\bigskip

\noindent (c) This follows from (b) and the proof of 
(a).\hfill{$\square$}\bigskip

 Now we come to some examples of converging situations, with some 
cycle relations involved.\bigskip  

\noindent {\bf Notation.}~Let $S\subseteq\omega$, and 
$$A^S:=\{ (s0\gamma ,s1\gamma )/~s\!\in\! 2^{<\omega}\hbox{\rm and\ 
Card}(s)\!\in\! S~\hbox{\rm and}~\gamma\!\in\! 2^\omega\}.$$
($\hbox{\rm Card}(s)$ is the number of ones in $s$.) We define 
partial homeomorphisms 
$${f^S_{n}:\bigcup_{s\in 2^n,~\hbox{\rm Card}(s)\in S}~N_{s0}\rightarrow
\bigcup_{s\in 2^n,~\hbox{\rm Card}(s)\in S}~N_{s1}}$$ 
by ${f^S_{n}(s0\gamma ):=s1\gamma}$. Notice that $A^S=A^{f^S}$ is 
Borel. One can show the existence of ${\mathfrak  A}:2^{\omega}\rightarrow 
2^\omega$ continuous such that ${\mathfrak  A}(S)$ is a Borel code for $A^S$, for 
each $S\subseteq\omega$. Notice that $\big(2^\omega ,(f^S_{n})_n\big)$ is 
a converging situation if and only if $S$ is infinite. This is also equivalent 
to $A^{S}\notin \mbox{pot}(\bormone )$. Indeed, if $S$ is finite, 
$\overline{A^{S}}\setminus A^{S}$ is a countable subset of $\Delta 
(2^\omega )$.\bf ~So in the sequel we will assume that $S$ is infinite.\rm
\bigskip

 Let $n_{S}:=\hbox{\rm min}~S$, and $S':=\{ n\! -\! n_{S}/n\!\in\! S\}$. Then 
$0\in S'$ and the maps $u$ and $v$ defined by $u(\alpha )=v(\alpha 
):=1^{n_{S}}\alpha$ are 
witnesses to $A^{S'}\leq^r_{c} A^{S}$.\bf ~So in the sequel we will 
also assume that $0\in S$.\rm\bigskip

\noindent $\bullet$ If $S\subseteq\omega$ and 
$t\in\omega^{<\omega}\setminus\{\emptyset\}$, then we set 
$f^S_{t}:=f^S_{t(0)}\ldots f^S_{t(|t|-1)}$, when it makes sense. We 
will also use the following tree 
${\mathfrak  R}$ on $2\times 2$. If $s,t\in 2^{<\omega}$, then we set
$$s~{\mathfrak  R}~t~\Leftrightarrow ~|s|\! =\! |t|~\mbox{and}~
(N_{s}\times N_{t})\cap A^S\!\not=\!\emptyset .$$
In particular, if $n_0\! <\! n_1$ and $1\!\in\! S$, then  we get 
${f^{S}_{<n_0,n_1>}(0^\infty )=f^{S}_{<n_1,n_0>}(0^\infty )}$.  
This is the kind of cycle relation we mentioned in the 
introduction. In this case $s({\mathfrak  R})$ is not acyclic since 
$<0^{n_1+1},0^{n_0}10^{n_1-n_0},0^{n_0}10^{n_1-n_0-1}1,0^{n_1}1,0^{n_1+1}>$ 
is an 
$s({\mathfrak  R})$-cycle. We set $f^{C}_{n}:={f^{S}_{n}}\vert_{C\cap 
{f^{S}_{n}}^{-1}(C)}$ for each Borel subset $C$ of $2^\omega$, when 
$S$ is fixed.

\vfill\eject

\noindent $\bullet$ Let $(H)$ be the following hypothesis on $S$:
$$(H)~\left\{\!\!\!\!\!\!
\begin{array}{ll}  
& ~~~~\mbox{Let}~C\in\borel (2^\omega )\mbox{,}~\sigma ~
\mbox{be~a~finer~topology~on}~C~\mbox{such~that}~\big([C,\sigma ],
(f^{C}_{n})_{n}\big)\cr    
& \mbox{is~a~converging~situation,}~l\mbox{, }p\in\omega\mbox{.~Then~we~can~find}~n\geq l~
\mbox{and}~\gamma\in D_{f^{C}_{n}}\cr    
& \mbox{with}~\mbox{Card}(\gamma\lceil n)+\big(S\cap [0,p]\big)=S\cap
\big(\mbox{Card}(\gamma\lceil n)+[0,p]\big). 
\end{array}
\right.$$
The next result will lead to a combinatorial condition on $S$ 
implying the minimality of $A^{S}$ among non-potentially closed sets.

\begin{thm} Let $S$ satisfy $(H)$, 
$B\!\in\!\borel (2^\omega )$, and $\tau$ a finer topology on $B$ such 
that $\big([B,\tau ],(f^{B}_{n})_{n}\big)$ is a converging situation. Then 
$A^{S}\leq^r_{c} A^{S}\cap B^2$.\end{thm}  

\noindent {\bf Proof.} Let $X:=[B,\tau ]$, $f_{n}:=f^{B}_{n}$. We are 
trying to build continuous maps $u,~v:2^\omega\rightarrow X$ such 
that $A^S=(u\times v)^{-1}(A^f)$. We will actually have more: $u=v$ 
will be one-to-one. We set $s\wedge t\!:=\! s\lceil 
\mbox{max}\{n\!\in\!\omega /s\lceil n\! =\! t\lceil n\}$, for $s,t\in 
2^{<\omega}$.\bigskip

\noindent $\bullet$ We construct a sequence $(U_{s})_{s\in 
2^{<\omega}}$ of nonempty clopen subsets of $X$, $\phi: 
\omega\rightarrow\omega$ strictly increasing, and $\theta: 
\omega\rightarrow\omega$ such that
$$\begin{array}{ll} 
(i) & ~~U_{s^\frown i }\subseteq U_{s}.\cr    
(ii) & ~\hbox{\rm diam}(U_{s})\leq 1/|s|~\mbox{if}~s\not=\emptyset .\cr    
(iii) & (s~{\mathfrak  R}~t~\hbox{\rm and}~s\!\not=\! t)\Rightarrow\left\{\!\!\!\!\!\!
\begin{array}{ll}    
& U_{t}\! =\! f_{\phi (| s\wedge t| )}[U_{s}]\mbox{,}\cr   
& \theta (| s\!\wedge\! t| )\! +\!\big(S\cap [0,| s\!\wedge\! t| 
]\big)\! =\! S\cap\big( 
\theta (| s\!\wedge\! t| )\! +\! [0,| s\!\wedge\! t| ]\big)\mbox{,}\cr   
& \forall z\!\in\! U_{s}~~\mbox{Card}\big(z\lceil\phi(| s\!\wedge\! t| 
)\big)\! 
=\!\theta (|s\!\wedge\! t|)\! +\!\mbox{Card}(s\lceil |s\!\wedge\! t|). 
\end{array}
\right.\cr   
(iv) & ~(\neg ~s~{\mathfrak R}~t~\hbox{\rm and}~|s|=|t|)~\Rightarrow 
~(U_{s}\times U_{t})\cap [\bigcup_{q<|s|} 
\hbox{Gr}(f_{q})\cup\Delta (X)]=\emptyset .
\end{array}$$
$\bullet$ First we show that this construction is sufficient to get 
the theorem. We define a continuous map ${u:2^\omega\rightarrow X}$ 
by  $\{ u(\alpha )\}:=\bigcap_{n}U_{\alpha\lceil n}$. If 
$\alpha <_{\mbox{lex}}\beta$, then we 
have $\neg\beta\lceil r~{\mathfrak  R}~\alpha\lceil r$ if $r$ is big 
enough, thus by 
condition (iv), $\big( u(\beta ), u(\alpha )\big)$ is in 
${U_{\beta\lceil r}
\times U_{\alpha\lceil r}\subseteq X^2\setminus\Delta (X)}$. 
Therefore $u$ is 
one-to-one. If ${(\alpha ,\beta )\in A^{S}}$, fix $n$ such that 
$\beta =f^{S}_{n}(\alpha )$. Then $\alpha\lceil r$ and $\beta\lceil 
r$ satisfy 
the hypothesis in condition (iii) for each $r>n$. Therefore 
$u(\beta )=f_{\phi (n)}\big(u(\alpha )\big)$ and 
$\big(u(\alpha ),u(\beta )\big)\in A^{f}$. If 
$\alpha =\beta$, then $(\alpha ,\beta )\notin A^{S}$ and  
${\big(u(\alpha ),u(\beta )\big)\in\Delta (X)\subseteq\neg A^{f}}$. 
Otherwise, 
$(\alpha ,\beta )\notin\overline{A^{S}}$ and there is $r_{0}$ such 
that 
$\alpha\lceil r$ and $\beta\lceil r$ satisfy the hypothesis in 
condition (iv) 
for $r\geq r_{0}$. This shows that $\big(u(\alpha ),u(\beta )\big)\notin A^{f}$. 
So it is enough to do the construction.\bigskip

\noindent $\bullet$ We set $U_{\emptyset}:=X$. Suppose that 
$(U_{s})_{s\in 2^{\leq p}}$, $\big(\phi (j)\big)_{j<p}$ and 
$\big(\theta (j)\big)_{j<p}$ satisfying conditions (i)-(iv) have 
been constructed, which is done for $p=0$.\bigskip 

\noindent $\bullet$ We will use the relation ${\mathfrak  R}_\Theta$ 
defined before Proposition 18 with $\theta_n\! :=\! 0^n$. Notice that 
${\mathfrak  R}_\Theta\subseteq {\mathfrak  R}$. We set 
$t_{0}\!:=\!\theta_p0$. We define a partition of $2^{p+1}$ as 
follows. Using Proposition 18.(b) we set, for $k\!\in\!\omega$,
$$H_k\!:=\!\{ t\!\in\! 2^{p+1}/|p_{t,t_0}|\! =\! k\! +\! 1\}.$$ 
If $H_{k+1}$ is non-empty, then $H_k$ is non-empty. Thus we can find an 
integer $q$ such 
that $H_0$, ..., $H_q$ are not empty and $H_k$ is empty if $k>q$. We 
order 
$2^{p+1}$ as follows: $t_{0}$, then $H_{1}$ in any order with 
$\theta_p1$ 
first, $H_{2}$ in any order, $\ldots$, $H_{q}$ in any order. This 
gives $t_{0}$, $\ldots$, $t_{2^{p+1}-1}$. Notice that we can find 
$j<n$ such that 
$t_{j}~s({\mathfrak  R}_\Theta )~t_{n}$ if $0<n<2^{p+1}$. In particular, 
if $E^n:=\{t_{j}/j\leq n\}$, then $\big( E^n,s({\mathfrak  R}_\Theta )\big)$ is 
connected for each $n<2^{p+1}$.

\vfill\eject 

\noindent $\bullet$ We will construct integers $\phi (p)$, $\theta 
(p)$ and nonempty 
clopen subsets $U^n_{k}$ of $X$, for $n<2^{p+1}$ and $k\leq n$, 
satisfying
$$\begin{array}{ll}  
(1) & U^n_{k}\subseteq U_{t_{k}\lceil p}.\cr   
(2) & \hbox{\rm diam}(U^n_{k})\leq 1/{p+1}.\cr   
(3) & (t_{k}~{\mathfrak  R}~t_{l}~\hbox{\rm and}~t_{k}\!\not=\! 
t_{l})\Rightarrow\cr  
& \ \ \ \ \ \ \ \ \ \ \ \ \ \ \ \ \ \ \ \ \left\{\!\!\!\!\!\!
\begin{array}{ll}   
& U^n_{l}\! =\! f_{\phi (|t_{k}\wedge t_{l}|)}[U^n_{k}]\mbox{,}\cr   
& \theta (|t_{k}\!\wedge\! t_{l}|)\! +\!\big(S\!\cap\! 
[0,|t_{k}\!\wedge\! t_{l}|]\big)\! =\! 
S\!\cap\!\big(\theta (|t_{k}\!\wedge\! t_{l}|)\! +\! 
[0,|t_{k}\!\wedge\! t_{l}|]\big)\mbox{,}\cr &
\forall z\!\in\! U^n_{k}~~\hbox{\rm Card}\big(z\lceil\phi(|t_{k}\!\wedge\! 
t_{l}|)\big)\! =\!  
\theta (|t_{k}\!\wedge\! t_{l}|)\! +\!\hbox{\rm  
Card}(t_{k}\lceil|t_{k}\!\wedge\! 
t_{l}|).
\end{array}\right.\cr   
(4) & \neg ~t_{k}~{\mathfrak  R}~t_{l}~\Rightarrow ~(U^n_{k}\times 
U^n_{l})\cap 
[\bigcup_{q\leq p} \hbox{Gr}(f_{q})\cup\Delta (X)]=\emptyset .\cr   
(5) & U^{n+1}_{k}\subseteq U^n_{k}. 
\end{array}$$
We will then set $U_{t_{k}}:=U^{2^{p+1}-1}_{k}$ for $k<2^{p+1}$, so 
that conditions (i)-(iv) are fullfilled.\bigskip

\noindent $\bullet$ Let $C\in\borone (U_{t_{0}\lceil 
p})\setminus\{\emptyset\}$ 
such that $C^{2}\cap\bigcup_{q\leq p}~\hbox{Gr}(f_{q})=\emptyset$. Apply 
hypothesis 
$(H)$ to $C$ and $\sigma\!:=\!\tau$. This gives 
${n_{0}\!\geq\!\hbox{\rm sup}\{\phi (q)\! +\! 1/q\! <\! p\}}$ and 
$\gamma\in D_{f^{C}_{n_{0}}}$ such that 
$$\hbox{\rm Card}(\gamma\lceil n_{0})+\big(S\cap [0,p]\big)=
S\cap\big( \hbox{\rm Card}(\gamma\lceil n_{0})+[0,p]\big).$$ 
We set $\phi (p):=n_{0}$, $\theta (p):=\hbox{\rm Card}(\gamma\lceil n_{0})$.\bigskip 

 We then choose $U^{0}_{0}\in\borone\big( C\cap 
f^{-1}_{n_{0}}(C)\big)\setminus\{\emptyset\}$ with suitable diameter 
such that $f_{n_{0}}[U^{0}_{0}]\cap U^{0}_{0}=\emptyset$, and 
$z\lceil n_{0}\! =\!\gamma\lceil n_{0}$ for each $z\!\in\! 
U^{0}_{0}$. Assume that ${U^{0}_{0}, ~\ldots,~U^{n-1}_{0},}$ 
${\ldots, U^{n-1}_{n-1}}$ satisfying conditions (1)-(5) have been 
constructed (which has already been accomplished for $n=1$). As $n\geq 1$, we have 
$t_{n}\not= t_{0}$ and $|p_{t_n,t_0}|\geq 2$. So fix $r<n$ such that 
$p_{t_n,t_0}(1)=t_{r}$. Notice that $U^{n-1}_{r}$ has been 
constructed.\bigskip

\noindent {\bf Case 1.}~$t_{n}\lceil p=t_{r}\lceil p$.\bigskip

\noindent - We have $t_{n}\lceil p=\theta_p$, thus $|p_{t_n,t_0}| 
=2$, $r=0$, 
$t_{n}=\theta_p1$ and $n=1$. Moreover, $U^{0}_{0}$ is a subset of 
${f^{-1}_{\phi (p)}}(U_{t_{1}\lceil p})$, so we can choose a nonempty 
clopen subset $U^1_{1}$ of $f_{\phi (p)}[U^{0}_{0}]$ with suitable 
diameter. 
Then we set $U^1_{0}\!:=\! f^{-1}_{\phi (p)}(U^1_{1})\!\subseteq\! 
U^{0}_{0}$. So 
${U^n_{0},\ldots ,U^n_{n}}$ are constructed and fullfill (1)-(3) and 
(5). It remains to check condition (4).\bigskip

\noindent - Fix ${k,~l\leq 1}$ such that $\neg ~t_{k}~{\mathfrak  
R}~t_{l}$. Then 
$k=1=1-l$. We have $U^1_{1} = f_{\phi(p)}[U^1_{0}]$. Thus 
$$U^1_{1}\times U^1_{0}=
f_{\phi(p)}[U^1_{0}]\times U^1_{0}=
f_{n_{0}}[U^1_{0}]\times U^1_{0}\subseteq 
f_{n_{0}}[U^0_{0}]\times U^0_{0}\subseteq C^{2}\mbox{,}$$ 
so we are done by the choice of $C$ and $U^{0}_{0}$.\bigskip

\noindent {\bf Case 2.}~$t_{n}\lceil p\not=t_{r}\lceil p$.\bigskip

\noindent 2.1. $t_{r}~{\mathfrak  R}_{\Theta}~t_{n}$.\bigskip

\noindent - By the induction hypothesis we have $U_{t_{n}\lceil p}=
f_{\phi(|t_{r}\wedge t_{n}|)}[U_{t_{r}\lceil p}]$ and 
$U^{n-1}_{r}\subseteq U_{t_{r}\lceil p}$. We choose a nonempty clopen 
subset $U^n_{n}$ of $f_{\phi(| t_{r}\wedge t_{n}| )}[U^{n-1}_{r}]$ with 
suitable diameter, so that conditions (1)-(5) for $k = l = n$ are 
fullfilled.\bigskip

\noindent - We then define the $U_{q}^n$'s for $q < n$, by induction 
on $|p_{t_{q},t_{n}}|$: fix $m\leq n$ with $p_{t_{q},t_{n}}(1) = t_{m}$. 
Notice that $q=r$ if $m=n$.

\vfill\eject

\noindent 2.1.1. $t_{m}~{\mathfrak R}_{\Theta}~t_{q}$.\bigskip

 We have $m<n$ since we cannot have $p_{t_{q},t_{n}}(1)~{\mathfrak 
R}_{\Theta}~t_{q}$  and 
$t_{q}~{\mathfrak R}_{\Theta}~p_{t_{q},t_{n}}(1)$ ($\tilde s \leq_{\hbox{\rm lex}} 
\tilde t$ if 
$\tilde s~{\mathfrak R}_{\Theta}~\tilde t$). So 
$U^{n-1}_{q} = f_{\phi (|t_{m}\wedge t_{q}|)}[U^{n-1}_{m}]$. We put
$$U_{q}^n:= f_{\phi (|t_{m}\wedge t_{q}|)}[U_{m}^n].$$
The set $U^{n}_{q}$ is a nonempty clopen subset of $U^{n-1}_{q}$ 
since $U^{n}_{m}\subseteq U^{n-1}_{m}$.\bigskip

\noindent 2.1.2. $t_{q}~{\mathfrak R}_{\Theta}~t_{m}$.\bigskip

 If $m<n$, then we have $U^{n-1}_{m} = f_{\phi (|t_{q}\wedge 
t_{m}|)}[U^{n-1}_{q}]$. We put
$$U_{q}^n:= f_{\phi (|t_{q}\wedge t_{m}|)}^{-1}(U_{m}^n)\mbox{,}$$
so that $U^{n}_{q}$ is a nonempty clopen subset of $U^{n-1}_{q}$. If 
$m=n$, then $q=r$ and the same conclusion holds, by the choice of 
$U^{n}_{n}$.\bigskip 

\noindent - So condition (5) is fullfilled in both cases. Conditions (1) and (2) 
are fullfilled for $k=q$, too. Let us check that the first part of condition (3) 
restricted to ${\mathfrak  R}_{\Theta}$ is fullfilled. Fix $k\not= l\leq n$ with 
$t_{k}~{\mathfrak R}_{\Theta}~t_{l}$. If $|p_{t_{k},t_{n}}| = 1$ 
and $|p_{t_{l},t_{n}}| = 2$, then the link between $t_{k}$ and $t_{l}$ has 
already been considered. The argument is similar if $|p_{t_{k},t_{n}}| = 2$ and 
$|p_{t_{l},t_{n}}| = 1$. If 
$|p_{t_{k},t_{n}}|$ and $|p_{t_{l},t_{n}}|$ are at least 2, then 
$p_{t_{k},t_{n}}(1) = p_{t_{l},t_{n}}(0)$ or 
$p_{t_{k},t_{n}}(0) = p_{t_{l},t_{n}}(1)$, by Proposition 18.(b). Here 
again, the link has already been considered. So 
condition (3) restricted to ${\mathfrak  R}_{\Theta}$ is fullfilled. It 
remains to check conditions (3) and (4).\bigskip

\noindent - Fix $k\not= l$ such that $t_{k}~{\mathfrak  R}~t_{l}$. 
Then $t_{k}$, $t_{l}$ differ at one coordinate only, and 
$t_{k}<_{\mbox{lex}}t_{l}$.\bigskip

\noindent {\bf Claim.}~\it Assume that $t_{k}$, $t_{l}$ differ at one 
coordinate only, and that $t_{k}<_{\mbox{lex}}t_{l}$. Then 
$$\mbox{Card}\big(z\lceil\phi(|t_{k}\wedge t_{l}|)\big)=\theta (|t_{k}\wedge 
t_{l}|)+\mbox{Card}(t_{k}\lceil|t_{k}\wedge t_{l}|)$$ 
for each $z\in U^n_{k}$.\bigskip\rm

 We can write
$$\begin{array}{ll} 
t_{k} &:=0^{n_{0}}10^{n_{1}}1\ldots 0^{n_{j-1}}10^{n_{j}}1
0^{n_{j+1}}1\ldots 0^{n_{q-1}}10^{n_{q}}\mbox{,}\cr\cr    
t_{l} &:=0^{n_{0}}10^{n_{1}}1\ldots 0^{n_{j-1}}10^{n'_{j}}10^m
10^{n_{j+1}}1\ldots 0^{n_{q-1}}10^{n_{q}}~~(n'_{j}+1+m=n_{j}). 
\end{array}$$
\indent By construction we have\bigskip

\leftline{$\begin{array}{ll} 
U^n_{k} & \!\!\! =\! f_{\phi (n_{0})} f_{\phi (\Sigma_{r\leq 
1}~(n_{r}+1)-1)}\ldots f_{\phi (\Sigma_{r\leq 
q-1}~(n_{r}+1)-1)}[U^n_{0}]\mbox{,}\cr\cr  
U^n_{l} & \!\!\! =\! f_{\phi (n_{0})}\!\ldots\! 
f_{\phi (\Sigma_{r\leq j-1}~(n_{r}+1)-1)}
f_{\phi (\Sigma_{r\leq j-1}~(n_{r}+1)+(n'_{j}+1)-1)}
f_{\phi (\Sigma_{r\leq j}~(n_{r}+1)-1)}\!\ldots 
\end{array}$} \smallskip

\rightline{$f_{\phi (\Sigma_{r\leq q-1}~(n_{r}+1)-1)}[U^n_{0}].$}

\vfill\eject

Notice that the length of $t_{k}\!\wedge\! t_{l}$ is equal to 
${\Sigma_{r\leq j-1}~(n_{r}\! +\! 1)\! +\! (n'_{j}\! +\! 1)\! -\! 
1}$. Set 
$$f:=f_{\phi (n_{0})}\!\ldots\! f_{\phi (\Sigma_{r\leq j-1}~(n_{r}+1)-1)}.$$ 
Then $U^n_{l}=ff_{\phi (|t_{k}\wedge t_{l}|)}f^{-1}(U^n_{k})$. Fix  
$k'\not= l'\leq n$ such that 
$$\begin{array}{ll} 
t_{k'} 
& := 0^{\Sigma_{r\leq j-1}~(n_{r}+1)+n_{j}}10^{n_{j+1}}1\ldots 0^{n_{q-1}}10^{n_{q}},\cr\cr 
t_{l'} 
& := 0^{\Sigma_{r\leq j-1}~(n_{r}+1)+n'_{j}}10^m10^{n_{j+1}}1\ldots 0^{n_{q-1}}10^{n_{q}}. 
\end{array}$$ 
Note that ${\hbox{\rm Card}\big(y\lceil\phi(|t_{k}\wedge t_{l}|)\big) =  
\theta (|t_{k}\wedge t_{l}|)+\hbox{\rm Card}(t_{k'}\lceil |t_{k}\wedge 
t_{l}|)}$, for each $y\in U^n_{k'}$, since $t_{k'}~{\mathfrak  R}_{\Theta}~t_{l'}$. 
But ${\hbox{\rm Card}\big(z\lceil\phi(|t_{k}\wedge t_{l}|)\big) = 
\hbox{\rm Card}\big(y\lceil\phi(|t_{k}\wedge t_{l}|)\big)+j}$, for each $z$ in 
$U^n_{k}=f[U^n_{k'}]$. As 
$$\hbox{\rm Card}(t_{k'}\lceil |t_{k}\wedge t_{l}|) = 
\hbox{\rm Card}(t_{k}\lceil |t_{k}\wedge t_{l}|)-j\mbox{,}$$ 
we get
$$\hbox{\rm Card}\big(z\lceil\phi(|t_{k}\wedge t_{l}|)\big)=\theta 
(|t_{k}\wedge t_{l}|)+\hbox{\rm Card}(t_{k}\lceil |t_{k}\wedge 
t_{l}|).\leqno{(+)}$$
This proves the claim.\hfill{$\diamond$}\bigskip 

\noindent - The second assertion in condition (3) is clearly 
fullfilled since 
$|t_{k}\wedge t_{l}|=|t_{k'}\wedge t_{l'}|$. As $t_{k}~{\mathfrak  
R}~t_{l}$ and 
$t_{k}\not= t_{l}$ we get $\hbox{\rm Card}(t_{k}\lceil| t_{k}\wedge t_{l}| )\in 
S$. This implies that $S$ contains 
${\theta (|t_{k}\wedge t_{l}|)+\hbox{\rm Card}(t_{k}\lceil |t_{k}\wedge 
t_{l}|)}$. By the claim we get
$$U^n_{l}=f_{\phi(|t_{k}\wedge t_{l}|)}ff^{-1}(U^n_{k})=
f_{\phi(|t_{k}\wedge t_{l}|)}[U^n_{k}]$$ 
(the compositions $ff_{\phi(|t_{k}\wedge t_{l}|)}f^{-1}$ and 
$f_{\phi(|t_{k}\wedge t_{l}|)}ff^{-1}$ are defined on $U^n_{k}$, so 
they are equal on this set). Thus condition (3) is fullfilled.\bigskip 

\noindent - To get condition (4), fix $k,l\leq n$ with $\neg 
~t_{k}~{\mathfrak R}~t_{l}$, 
$v(i):=|p_{t_k,t_l}(i)\wedge p_{t_k,t_l}(i+1)|$, and $\varepsilon 
(i):=1$ (resp., $-1$) if 
$p_{t_k,t_l}(i)~{\mathfrak  R}_{\Theta}~p_{t_k,t_l}(i+1)$ (resp., 
$p_{t_k,t_l}(i+1)~{\mathfrak  R}_{\Theta}~p_{t_k,t_l}(i)$), for 
$i+1<|p_{t_k,t_l}|$. We set ${f^\varepsilon_{v}:={f_{\phi 
(v(|v|-1))}^{\varepsilon (|v|-1)}\ldots 
f_{\phi (v(0))}^{\varepsilon (0)}}}$, so that 
$U^n_{l}=f^\varepsilon_{v}(U^n_{k})$. 
Let $m$ be maximal such that $t_{k}(m)\not= t_{l}(m)$. As $\phi$ is 
strictly 
increasing, we get $(U^n_{k}\times U^n_{l})\cap\Delta (X)=\emptyset$, 
by Proposition 18.(c).\bigskip 

\noindent - If $t_{k}$, $t_{l}$ differ in at least two coordinates $m\not= m'$, 
then the number of appearances of $m$ and $m'$ in $v$ is odd. As $\phi$ is 
strictly increasing, this is also true for $\phi (m)\not=\phi (m')$ in 
$\{\phi\big( v(i)\big)/i<|v|\}$. This implies that 
$(U^n_{k}\times U^n_{l})\cap [\bigcup_{q\leq p}~\hbox{Gr}(f_{q})]=\emptyset$.\bigskip 

\noindent - If $t_{k}$, $t_{l}$ differ at only one coordinate $m$ and 
$t_{k}>_{\hbox{\rm lex}}t_{l}$, then ${\alpha\big(\phi 
(m)\big)>\beta\big(\phi (m)\big)}$ 
if $(\alpha ,\beta )$ is in $U^n_{k}\times U^n_{l}$, and 
$(U^n_{k}\times U^n_{l})\cap [\bigcup_{q\leq 
p}~\hbox{Gr}(f_{q})]=\emptyset$.\bigskip

\noindent - So we may assume that $t_{k}$, $t_{l}$ differ only at 
coordinate $\phi^{-1}(q)$, and that $t_{k}<_{\hbox{\rm lex}}t_{l}$. By the Claim we 
have $(+)$ for 
each $z\in U^n_{k}$. But $\hbox{\rm Card}(t_{k}\lceil |t_{k}\wedge t_{l}|)\notin 
S$, since 
$\neg ~t_{k}~{\mathfrak  R}~t_{l}$. So $\hbox{\rm Card}(z\lceil q)\notin S$ if $z\in 
U^n_{k}$, and $f_{q}$ is not defined on $U^n_{k}$.

\vfill\eject

\noindent 2.2. $t_{n}~{\mathfrak  R}_{\Theta}~t_{r}$.\bigskip

 This cannot hold since $t_{r}{\mathfrak  R}_{\Theta}t_{n}$. Indeed, if 
$t_{n}\! =\! 0^{n_{0}}10^{n_{1}}1\ldots 0^{n_{q-1}}10^{n_{q}}$, then
$$\begin{array}{ll} 
p_{t_n,t_0}(1) 
& = 0^{n_{0}+n_{1}+1}1\ldots 0^{n_{q-1}}10^{n_{q}}\mbox{,}\cr 
& \ \ .\cr 
& \ \ .\cr 
& \ \ .\cr 
p_{t_n,t_0}(|p_{t_n,t_0}|-2)\!\!\!\! 
& = 0^{n_{0}+n_{1}+\ldots +n_{q-1}+q-1}10^{n_{q}}. 
\end{array}$$
This finishes the proof.\hfill{$\square$}

\begin{lem} The set $S$ satisfies hypothesis $(H)$ if the 
following is fullfilled:
$$\forall p\!\in\!\omega~\exists k\!\in\!\omega~\forall 
q\!\in\!\omega~\exists c\!\in\!\omega\!\cap\! 
[q,q\! +\! k]~~c\! +\!\big(S\!\cap\! [0,p]\big)\! =\! S\!\cap\!\big( 
c\! +\! [0,p]\big).\leqno{(M)}$$
In particular, condition (M) implies that $A^{S}$ is minimal among 
non-potentially closed sets for $\leq^r_{c}$.\end{lem} 

\noindent {\bf Proof.} Note that $\emptyset\not=\Delta (C)
\subseteq\overline{\bigcup_{q\geq l}~\hbox{Gr}(f^{C}_{q})}$, since 
$\big([C,\sigma ],
(f^{C}_{n})_{n}\big)$ is a converging situation. So fix $q_{0}\geq 
l$ such that 
$D_{f^{C}_{q_{0}}}\not=\emptyset$, and $O_{0}:=D_{f^{C}_{q_{0}}}$. 
Assume that $q_{r}$ and $O_{r}$ have been constructed. We then choose 
$q_{r+1}>q_{r}$ such that ${O_{r}\cap 
(f^{C}_{q_{r+1}})^{-1}(O_{r})\not=\emptyset}$,  
and we define ${O_{r+1}:=O_{r}\cap (f^{C}_{q_{r+1}})^{-1}(O_{r})}$. 
This gives $(q_{r})_{r<M}$ and $(O_{r})_{r<M}$, where 
$M:=p+k$.\bigskip

\noindent $\bullet$ For $t\in\omega^{<\omega}$, we let 
$f^C_{t}:=f^C_{t(0)}\ldots f^C_{t(|t|-1)}$, when it makes sense. We choose 
$${n\geq\hbox{\rm max}(\hbox{\rm sup}\{q_{r}+1/r<M\},l)}$$ 
with ${f^{C}_{q_{0},\ldots ,q_{M-1}}[O_{M-1}]\cap {f^{C}_{n}}^{-1}(f^{C}_{q_{0},\ldots 
,q_{M-1}}[O_{M-1}])}\not=\emptyset$. Let $\beta:=f^{C}_{q_{0},\ldots ,q_{M-1}}(\alpha )$ be in the  
intersection. Notice that $q:=\hbox{\rm Card}(\beta\lceil n)-M\in\omega$. 
This gives $c$ in $\omega\cap [q,q+k]$. As $0\in S$, there is $j\leq k$  
with ${c=\hbox{\rm Card}(\beta\lceil n)-p-j\in S}$. Notice that 
$\beta = f^{C}_{q_{0},\ldots ,q_{p+j-1}}(\gamma )$, where 
$\gamma = f^{C}_{q_{p+j},\ldots ,q_{M-1}}(\alpha )$. As 
${\hbox{\rm Card}(\gamma\lceil n) =c}$, $f^{S}_{n}(\gamma )$ 
is defined. But $f^{C}_{n}(\beta )$ is in $f^{C}_{q_{0},\ldots 
,q_{M-1}}[O_{M-1}]$ 
and $f^{S}_{n}(\gamma )$ is in $C$. So $f^{C}_{n}(\gamma )$ is 
defined.\bigskip

\noindent $\bullet$ The lemma now follows from Corollary 16 and Theorem 
19.\hfill{$\square$}\bigskip  

\noindent {\bf Example.} ~We set ${S_{m,F}:=\{n\!\in\!\omega / 
n~(\hbox{\rm mod}~m)\!\in\!\{ 0\}\!\cup\! F\}}$, where 
$m\in\omega\setminus\{ 0\}$ and ${F\subseteq m\setminus\{ 0\}}$. Then 
$S_{m,F}$ fullfills condition (M). In particular, $A^{\omega}$ is 
minimal. But this gives only countably many examples. To get more, we 
need some more notation:\bigskip  

\noindent {\bf Notation.} ~For $\beta\in \omega^\omega$, we set 
$S_{\beta}:=\{\Sigma_{i<l}~\big(1+\beta (i)\big)/l\!\in\!\omega\}$. 
Notice that $0\in S_{\beta}$, $S_{\beta}$ is infinite, and that any 
infinite $S$ containing $0$ is of this form. Moreover, the map 
$\beta\mapsto S_{\beta}$ is continuous since 
$n\!\in\! S_{\beta}\Leftrightarrow\exists l\!\leq\! n~~n\! =\! 
\Sigma_{i<l}~\big(1\! +\!\beta (i)\big)$. We will define a family 
$(\beta_{\alpha})_{\alpha\in 2^\omega}$. Actually, we can find at 
least two examples:\bigskip

\noindent $\bullet$ The original example is the following. For $\alpha\in 
2^\omega$, we recursively define a sequence 
$(s_{\alpha ,n})_{n}\subseteq 2^{<\omega}$ as 
follows: $s_{\alpha ,0}:=0$, $s_{\alpha ,1}:=1$, $s_{\alpha ,n+2}:=
s_{\alpha ,n}^{\alpha (n)+1}s_{\alpha ,n+1}^{\alpha (n+1)+1}$. Notice 
that $s_{\alpha ,n}\prec_{\not=}s_{\alpha ,n+2}$, so that 
$\beta_{\alpha}:= \hbox{\rm lim}_{n\rightarrow\infty}~s_{\alpha ,2n}\in 
2^\omega$ is defined.\bigskip

\noindent $\bullet$ A. Louveau found another example for which it is simpler to 
check property (M) (and ($\perp$) later), and in the sequel we will work 
with it. For $\alpha\in 2^\omega$, $n\in\omega$ and $\varepsilon\in 2$, we set 
$$\gamma_{\alpha}(4n+2\varepsilon ):=\varepsilon\mbox{,}$$ 
$\gamma_{\alpha}(2n+1):=\alpha (n)$ (so that $\gamma_{\alpha}$ has infinitely many zeros and ones, and the map $\alpha\mapsto\gamma_{\alpha}$ is continuous). For $i\in\omega$, 
we then set $(i)_{0}:=\hbox{\rm max}\{m\!\in\!\omega /2^m~\hbox{\rm divides}~i+1\}$. 
Finally, we put $\beta_{\alpha}(i):=\gamma_{\alpha}\big((i)_{0}\big)$.\bigskip 

 Notice that the map $\alpha\mapsto\beta_{\alpha}$ is continuous, so that the 
map $\alpha\mapsto A^{S_{\beta_{\alpha}}}$ is continuous in the codes.

\begin{cor}  Let $\alpha\in 2^\omega$. Then 
$S_{\beta_{\alpha}}$ satisfies condition (M). In particular, 
$A^{S_{\beta_{\alpha}}}$ is minimal among non-potentially closed sets 
for $\leq^r_{c}$.\end{cor}

\noindent {\bf Proof.} First notice that it is 
enough to show that the following is fullfilled:
$$\forall P\!\in\!\omega~\exists K\!\in\!\omega~
\forall Q\!\in\!\omega~\exists C\!\in\!\omega\!\cap\! [Q,Q\! +\! K]~~
\beta_{\alpha}\lceil P\!\prec\!\beta_{\alpha}\! -\! 
\beta_{\alpha}\lceil C.\leqno{(MM)}$$
Indeed, this condition associates $K$ to $P:=p$. Set 
$k:=2K+1$. For $q\in\omega$, let $Q$ be minimal with 
$\Sigma_{i<Q}~\big(1+\beta_{\alpha} (i)\big)\geq q$, and fix 
$C\in\omega\cap 
[Q,Q+K]$ such that $\beta_{\alpha}\lceil P\prec\beta_{\alpha}-
\beta_{\alpha}\lceil C$. We put 
$$c:=\Sigma_{i<C}~\big(1+\beta_{\alpha} (i)\big).$$ 
Notice that $c\leq q+k$ since
$c\!\leq\!\Sigma_{i<Q-1}~\big(1\! +\!\beta_{\alpha} (i)\big)
\! +\!\Sigma_{Q-1\leq i<C}~\big(1\! +\!\beta_{\alpha} (i)\big)\! <\! 
q\! +\! 2(C\! -\! Q\! +\! 1)\leq q\! +\! 2(K\! +\! 1)$. 
Finally, note that 
$c+\Sigma_{i<j}~\big(1+\beta_{\alpha} (i)\big)=\Sigma_{i<C+j}~\big(1+\beta_{\alpha} (i)\big)$, by induction on $j\leq p$.\bigskip

 Notice that for any integers $n,~i$ and $l$ with $i<2^n-1$, we have 
$(2^n\cdot l+i)_{0}=(i)_{0}$. Indeed, we can find $N$ with 
$i=2^{(i)_{0}}(2N+1)-1$, 
and $(i)_{0}<n$. Thus ${2^n\cdot l+i=2^{(i)_{0}}(2^{n-(i)_{0}}\cdot l+2N+1)-1}$ 
and $(2^n\cdot l+i)_{0}=(i)_{0}$. Now, if $P\in\omega$, then let $n_{0}$ be 
minimal 
with $K:=2^{n_{0}}-1\geq P$. If $Q\in\omega$, then let 
$l\in\omega\cap [\frac{Q}{2^{n_{0}}},\frac{Q}{2^{n_{0}}}+1[$ and 
$C:=2^{n_{0}}\cdot l$. If $i<P$, then $i<2^{n_{0}}-1$, so 
${(2^{n_{0}}\cdot l+i)_{0}=(i)_{0}}=(C+i)_{0}$. Thus 
$\beta_{\alpha}(i)=\beta_{\alpha}(C+i)$.\hfill{$\square$}\bigskip

Now we come to the study of the cardinality of complete families of minimal sets.

\begin{lem} Let $\big(X,(f_{n})\big),\big(X',(f'_{n})\big)$ be 
converging situations, and ${u,v:X\!\rightarrow\! X'}$ continuous 
maps such that $A^{f}=(u\times v)^{-1}(A^{f'})$. Then 
$u=v$.\end{lem}

\noindent {\bf Proof.} For $x\in X$, fix $x_{k}\in X$ and 
$n_{k}\in\omega$ such 
that $\big(x_{k},f_{n_{k}}(x_{k})\big)$ tends to $(x,x)$. Note that 
$\big(u(x),v(x)\big)\notin A^{f'}$. Moreover,  
$\big(u[x_{k}],v[f_{n_{k}}(x_{k})]\big)\in A^{f'}$. Thus 
${\big(u(x),v(x)\big)\in\overline{A^{f'}}\setminus A^{f'}=\Delta 
(X')}$, therefore $u=v$.\hfill{$\square$}\bigskip

 Recall that $A^{-1}:=\{(y,x)\in X^2/(x,y)\in A\}$ if $A\subseteq 
X^2$.

\begin{thm} Fix $S,~S'$ satisfying condition (M). Then\smallskip  

\noindent (a) $A^{S}\perp^r_{c}A^{S'}$, provided that the following condition is fullfilled:
$$\exists p\!\in\!\omega ~\forall c\!\in\!\omega ~~c\! 
+\!\big(S\!\cap\! [0,p]\big)\!\not=\! S'\!\cap\!\big(c\! +\! 
[0,p]\big).\leqno{(\perp )}$$
(b) $A^{S}\perp^r_{c}(A^{S'})^{-1}$, provided that the following 
condition is fullfilled:
$$\exists p\!\in\!\omega ~\forall c\!\in\!\omega ~~c\! 
-\!\big(S\!\cap\! [0,p]\big)\!\not=\! S'\!\cap\!\big(c\! -\! 
[0,p]\big).\leqno{(\perp^{-1})}$$\end{thm} 

\noindent {\bf Proof.} (a) We argue by contradiction: by Lemma 20, we 
can find continuous maps ${u,~v:2^\omega\!\rightarrow\! 2^\omega}$ 
such that $A^{S}=(u\times v)^{-1}(A^{S'})$. By Lemma 22, we have 
$u=v$.\bigskip 

\noindent\bf Claim.\it\ Let $n$, $q$ be integers and 
$N\in\boraone(2^\omega)\setminus\{\emptyset\}$. Then we can find 
integers ${n'> n}$, ${q'\! >\! q}$ and a nonempty open subset $N'$ of 
$N\cap {f^S_{n'}}^{-1}(N)$ with 
${{f^{S'}_{q'}}[u(\gamma )]\! =\! u[f^S_{n'}(\gamma )]}$, and 
$${\mbox{Card}(\gamma\lceil n')+\big(S\cap [0,p]\big)=
S\cap\big(\mbox{Card}(\gamma\lceil n')+[0,p]\big),}$$
for each $\gamma\in N'$.\bigskip\rm

 Indeed, let ${\delta\in u[N]}$. As $(\delta ,\delta )$ is not in 
$\bigcup_{q'\leq q}\hbox{Gr}(f^{S'}_{q'})$, we can find a clopen 
neighborhood $W$ of $\delta$ such that ${W^2\cap\bigcup_{q'\leq q} 
\hbox{Gr}(f^{S'}_{q'})=\emptyset}$. 
Let $\tilde N\in\borone (2^\omega )\setminus\{\emptyset\}$ with 
${\tilde N\subseteq N\cap u^{-1}(W)}$. By Lemma 20, we can find 
$n'\! >\! n$ and $\gamma_{0}\!\in\!\tilde N\cap 
{f^S_{n'}}^{-1}(\tilde N)$ with
$${\hbox{\rm Card}(\gamma_{0}\lceil n')\! +\!\big(S\!\cap\! [0,p]\big)\! =\! 
S\!\cap\!\big(\hbox{\rm Card}(\gamma_{0}\lceil n')\! +\! [0,p]\big).}$$ 
Now there is $q'(\gamma )$ such that ${{f^{S'}_{q'(\gamma 
)}}[u(\gamma )] = 
u[f^S_{n'}(\gamma )]}$, for ${\gamma\in {\tilde N}\cap 
{f^S_{n'}}^{-1}({\tilde N})\cap N_{\gamma_{0}\lceil n'}}$. We have 
$q'(\gamma )>q$, by the choice of 
$W$. By Baire's Theorem we get $q'$ and 
$N'$.\hfill{$\diamond$}\bigskip 

 By the Claim we get $n_1$, $q_{1}$ and $N_{1}\subseteq 
D_{f^S_{n_1}}$ with 
${f^{S'}_{q_{1}}}[u(\gamma )] = u[f^S_{n_{1}}(\gamma )]$ and
$${\hbox{\rm Card}(\gamma\lceil n_{1})+\big(S\cap [0,p]\big)=
S\cap\big(\hbox{\rm Card}(\gamma\lceil n_{1})+[0,p]\big)\mbox{,}}$$ 
for each ${\gamma\in N_{1}}$.\bigskip

 We then get $n_2>n_{1}$, ${q_{2}>q_1}$, and a nonempty 
open subset $N_{2}$ of $N_{1}\cap {f^S_{n_{2}}}^{-1}(N_{1})$ with 
${{f^{S'}_{q_{2}}}[u(\gamma )] = u[f^S_{n_{2}}(\gamma )]}$ and
$${\hbox{\rm Card}(\gamma\lceil n_{2})+\big(S\cap [0,p]\big)=
S\cap\big(\hbox{\rm Card}(\gamma\lceil n_{2})+[0,p]\big)\mbox{,}}$$
for each $\gamma$ in $N_{2}$. We continue in this fashion, until we get 
$n_{p+1}$, $q_{p+1}$ and $N_{p+1}$. Fix  
$\gamma\!\in\! N_{p+1}$ and set $c:=\hbox{\rm Card}\big(u(\gamma )\lceil 
q_{p+1}\big)$.\bigskip

\noindent $\bullet$ Fix ${m\in S\!\cap\! [0,p]}$. For  
$t\in\omega^{<\omega}$, we set $f^S_{t}:=f^S_{t(0)}\ldots 
f^S_{t(|t|-1)}$, when it makes sense. Notice that 
${f^S_{n_{p-m+1},\ldots ,n_{p+1}}(\gamma ) = 
f^{S}_{n_{p+1},n_{p-m+1},\ldots ,n_{p}}(\gamma )}$ is defined. Therefore, $A^{S}$ contains 
$$\big(f^{S}_{n_{p-m+1},\ldots ,n_{p}}(\gamma ), f^{S}_{n_{p-m+1},\ldots ,n_{p+1}}(\gamma )\big)
\mbox{,}$$
which implies that $A^{S'}$ contains ${\big(u[f^{S}_{n_{p-m+1},\ldots ,n_{p}}(\gamma )], 
u[f^{S}_{n_{p-m+1},\ldots ,n_{p+1}}(\gamma )]\big)}$. This 
shows that $A^{S'}$ contains $\big(f^{S'}_{q_{p-m+1},\ldots 
,q_{p}}[u(\gamma )],
f^{S'}_{q_{p-m+1},\ldots ,q_{p+1}}[u(\gamma )]\big)$, thus 
$$f^{S'}_{q_{p-m+1},\ldots ,q_{p+1}}[u(\gamma )]\! =\! 
f^{S'}_{q_{p+1},q_{p-m+1},\ldots ,q_{p}}[u(\gamma )]\mbox{,}$$
so ${c+\big( S\cap [0,p]\big)\subseteq 
S'\cap\big(c+[0,p]\big)}$.\bigskip 

\noindent $\bullet$ Conversely, let $m:=c+m'\in 
S'\cap\big(c+[0,p]\big)$. Again 
$f^S_{n_{p-m'+1},\ldots ,n_{p+1}}(\gamma )$ is defined. Notice that
$$u[f^S_{n_{p-m'+1},\ldots ,n_{p+1}}(\gamma )]\! =\! 
f^{S'}_{q_{p-m'+1},\ldots ,q_{p+1}}[u(\gamma )]\! =\! 
f^{S'}_{q_{p+1},q_{p-m'+1},\ldots ,q_{p}}[u(\gamma )].$$ 
Therefore $\big(u[f^S_{n_{p-m'+1},\ldots ,n_{p}}(\gamma )],
u[f^S_{n_{p-m'+1},\ldots ,n_{p+1}}(\gamma )]\big)\!\in\! A^{S'}$, 
$A^{S}$ contains the pair 
$$\big(f^S_{n_{p-m'+1},\ldots ,n_{p}}(\gamma ),f^S_{n_{p-m'+1},\ldots ,n_{p+1}}(\gamma )\big)\mbox{,}$$ and $f^S_{n_{p-m'+1},\ldots ,n_{p+1}}(\gamma ) = f^S_{n_{p+1},n_{p-m'+1},\ldots ,n_{p}}(\gamma )$.  
Therefore $\hbox{\rm Card}(\gamma\lceil n_{p+1})+m'\!\in\! S$ and 
$m'\!\in\! S$, so ${S'\cap\big(c\! +\! [0,p]\big)\subseteq c\! +\!\big( S\cap [0,p]\big)}$. But this contradicts condition {($\perp$)} since we actually have the equality.\bigskip

\noindent (b) The proof is similar to that of (a). This 
time $A^{S}=(u\times v)^{-1}\big( (A^{S'})^{-1}\big)$. We construct 
sequences $(n_{j})_{1\leq j\leq p+1}$, $(q_{j})_{1\leq j\leq p+1}$ and 
$(N_{j})_{1\leq j\leq p+1}$ satisfying the equality 
${{f^{S'}_{q_{j}}}^{-1}[u(\gamma )] = u[f^S_{n_{j}}(\gamma )]}$ and
$${\hbox{\rm Card}(\gamma\lceil n_{j})+\big(S\cap 
[0,p]\big)=S\cap\big(\hbox{\rm Card}(\gamma\lceil n_{j})+[0,p]\big)\mbox{,}}$$ 
for each ${\gamma\in N_{j}}$. This gives
$$(f^{S'}_{q_{p-m+1}})^{-1}\ldots 
(f^{S'}_{q_{p+1}})^{-1}[u(\gamma )] =  
(f^{S'}_{q_{p+1}})^{-1}
(f^{S'}_{q_{p-m+1}})^{-1}\ldots 
(f^{S'}_{q_{p}})^{-1}[u(\gamma )]\mbox{,}$$
thus $c\! -\!\big( S\cap [0,p]\big)\subseteq S'\cap\big(c\! -\! [0,p]
\big)$, and we complete the proof as we did for (a).\hfill{$\square$}

\begin{cor}  Let $\alpha\not= \alpha'\in 
2^\omega$. Then $S_{\beta_{\alpha}}$, $S_{\beta_{\alpha'}}$ satisfy 
conditions (M), ($\perp$) and ($\perp^{-1}$). In particular, 
$A^{S_{\beta_{\alpha}}}\perp^r_{c}A^{S_{\beta_{\alpha'}}}$ 
and 
$A^{S_{\beta_{\alpha}}}\perp^r_{c}(A^{S_{\beta_{\alpha'}}})^{-1}$.\end{cor}

Theorem 5 is a corollary of this result. We saw 
that the map $\alpha\mapsto A^{S_{\beta_{\alpha}}}$ is 
continuous in the codes, and it is injective by Corollary 24. This implies that 
$(A^{S_{\beta_{\alpha}}})_{\alpha\in 2^\omega}$ is a perfect 
antichain for $\leq^r_{c}$ made of minimal sets (we use Corollaries 21 and 24).\bigskip

\noindent {\bf Proof.} If $s\in 2^{<\omega}$ and $t\in 
2^{\leq\omega}$, we say 
that $s\subseteq t$ if we can find an integer $l\leq |t|$ such that 
$s\prec t-t\lceil l$. We define $s^{-1}\in 2^{| s| }$ by  
$s^{-1}(i):=s(| s| -1-i)$, for $i<| s| $. We say that $s$ is $symmetric$ if $s=s^{-1}$.\bigskip 

\noindent $\bullet$ It is enough to prove the following condition:
$$\exists P\!\in\!\omega~\ \beta_{\alpha}\lceil 
P\!\not\subseteq\!\beta_{\alpha'}\ \ \hbox{\rm and}\ \ 
(\beta_{\alpha}\lceil 
P)^{-1}\!\not\subseteq\!\beta_{\alpha'}.\leqno{(\perp\perp)}$$  
Indeed, we will see that $(\perp\perp )$ implies $(\perp )$ and 
$(\perp^{-1})$ of 
Theorem 23. Condition $(\perp\perp )$ gives $P>0$. Let $p:=2P$ and 
$c\in\omega$. We argue by contradiction.\bigskip 

\noindent $(\perp )$ Assume that 
${c+\big(S_{\beta_{\alpha}}\cap [0,p]\big)=S_{\beta_{\alpha'}}\cap\big(c+[0,p]\big)}$. As 
$0\in S_{\beta_{\alpha}}$, we can find $l$ with 
$$c\! =\!\Sigma_{i<l}~\big(1\! +\!\beta_{\alpha'}(i)\big).$$ 
It is enough to prove that if $n\! <\! P$, then $\beta_{\alpha}(n)\! =\!\beta_{\alpha'}(l\! +\! n)$. 

\vfill\eject

We argue by induction on $n$.\bigskip

\noindent - Notice that $\beta_{\alpha}(0)=0$ is equivalent to $1\in 
S_{\beta_{\alpha}}$ and to 
$\beta_{\alpha'}(l)=0$. Therefore 
${\beta_{\alpha}(0)=\beta_{\alpha'}(l)}$.\bigskip

\noindent - Now suppose that $n+1<P$ and 
$\beta_{\alpha}(m)=\beta_{\alpha'}(l+m)$,  
for each $m\leq n$. As 
$$2+\Sigma_{m\leq n}~\big(1+\beta_{\alpha}(m)\big)\leq p\mbox{,}$$ 
we get $\beta_{\alpha}(n+1)=\beta_{\alpha'}(l+n+1)$.\bigskip 

\noindent $(\perp^{-1})$ Assume that ${c-\big(S_{\beta_{\alpha}}\cap 
[0,p]\big)=
S_{\beta_{\alpha'}}\cap\big(c-[0,p]\big)}$. Let $l':=l-P$ (as 
$2P-1$ or $2P$ is in $S_{\beta_{\alpha}}\cap [0,p]$, $c>2P-2$ and 
$l'\geq 
0$). As $(\beta_{\alpha}\lceil P)^{-1}\not\subseteq\beta_{\alpha'}$ 
we 
can find $n<P$ such that 
$\beta_{\alpha}(n)\not=\beta_{\alpha'}(l-1-n)$, since 
$(\beta_{\alpha}\lceil P)^{-1}\not\prec\beta_{\alpha'} 
-\beta_{\alpha'}\lceil l'$. 
We conclude as in the case $(\perp )$.\bigskip 

\noindent $\bullet$ First notice that $\beta_{\alpha}\lceil 
(2^n-1)=[\beta_{\alpha}
\lceil (2^n-1)]^{-1}$ for each integer $n$. Indeed, let $i<2^n-1$. It 
is enough to see that $(i)_{0}=(2^n-2-i)_{0}$. But we have 
$${2^n-2-i=2^n-2-2^{(i)_{0}}(2N+1)+1=2^{(i)_{0}}(2^{n-(i)_{0}}-2N-1)-1\mbox{,}}$$ 
so we are done, since $2^{n-(i)_{0}}-2N-1$ is odd and positive. So it 
is enough to find $n$ such that $\beta_{\alpha}\lceil 
(2^n-1)\not\subseteq\beta_{\alpha'}$.\bigskip 

\noindent $\bullet$ Let $n_{0}$ minimal with 
$\gamma_{\alpha} (n_{0})\!\not=\!\gamma_{\alpha'} (n_{0})$, and 
$n_{1}\! >\! n_{0}\! +\! 1$ with 
${\gamma_{\alpha'} (n_{0}+1)\!\not=\!\gamma_{\alpha'} (n_{1})}$. We 
put $n\!:=\! n_{1} + 2$. We argue by contradiction: we get $l$ with  
${\gamma_{\alpha}\big((i)_{0})\! =\!\gamma_{\alpha'}\big((l+i)_{0})}$, for each 
$i\! <\! 2^n-1$.\bigskip 

\noindent $\bullet$ Notice that for each $m<n-1$ we can find 
$i<2^{n-1}$ with $(l+i)_{0}=m$. Indeed, let 
$${N\in\omega\cap [2^{-m-1}(l+1)-2^{-1},2^{-m-1}(2^{n-1}+l+1)-2^{-1}[}.$$ 
It is clear that $i:=2^m(2N+1)-l-1$ is suitable.\bigskip 

\noindent $\bullet$ Let $M\geq n_{0}$ and $(\varepsilon_{j})_{j\leq 
M}\subseteq 2$ with 
$l=\Sigma_{j\leq M}~\varepsilon_{j}\cdot 2^j$. For $k\leq n_{0}$ we define 
$${i_{k}\!:=\!\Sigma_{j<k}~(1\! -\!\varepsilon_{j})\cdot 2^j +\varepsilon_{k}\cdot 2^k}.$$ 
Note that $i_{k}\! <\! 2^{k+1}$ and ${l\! +\! i_{k}\equiv 2^k\! -\! 1~(\hbox{\rm mod}~2^{k+1})}$. We will 
show the following, by induction on $k$:\bigskip

\noindent - The sequence $\big(\beta_{\alpha}(i)\big)_{i<2^{n-1},i\equiv 
2^k-1~(\hbox{\rm mod}~2^{k+1})}$ is constant with value 
$\gamma_{\alpha} (k)$, and equal to 
$${\big(\beta_{\alpha'}(l + i)\big)_{i<2^{n-1},i\equiv 
2^k-1~(\hbox{\rm mod}~2^{k+1})}}.$$
- The sequence $\big(\beta_{\alpha'}(l+i)\big)_{i<2^{n-1},i\equiv 
i_{k}~(\hbox{\rm mod}~2^{k+1})}$ is constant with value 
$\gamma_{\alpha'} (k)$.\bigskip

\noindent - The sequence $\big(\beta_{\alpha'}(l+i)\big)_{i<2^{n-1},i\equiv 
i_{k}+2^k~(\hbox{\rm mod}~2^{k+1})}$ is not constant.\bigskip

\noindent - $\varepsilon_{k}=0$ and $\gamma_{\alpha} (k)=\gamma_{\alpha'} 
(k)$.\bigskip

 This will give the desired contradiction with $k=n_{0}$. 

\vfill\eject

 So assume that these facts have been shown for $j<k\leq n_{0}$.\bigskip 

\noindent - The first point is clear.\bigskip 

\noindent - The second one comes from the fact that $l+i$ is of the form 
$2^k(2K+1)-1$ if 
$i\equiv i_{k}~(\hbox{\rm mod}~2^{k+1})$, since ${l+i_{k}\equiv 
2^k-1~(\hbox{\rm mod}~2^{k+1})}$.\bigskip 

\noindent - To see the third one, choose $i<2^{n-1}$ such that 
$(l+i)_{0}=n_{0}+1$ (or 
$n_{1}$). We have to see that $i\equiv i_{k}+2^k~(\hbox{\rm 
mod}~2^{k+1})$. We can find 
$(\eta_{j})_{j<n-1}$ with $i=\Sigma_{j<n-1}~\eta_{j}\cdot 2^j$, so that
$$l+i+1\equiv 
1+\Sigma_{j<k}~\eta_{j}\cdot 2^j+(\varepsilon_{k}+\eta_{k})\cdot 2^k~
(\hbox{\rm mod}~2^{k+1})\mbox{,}$$ 
by the induction hypothesis. This inductively shows that $\eta_{j}=1$ 
if $j<k$ and ${\eta_{k}=1-\varepsilon_{k}}$. Thus $i\equiv 
2^k-1+(1-\varepsilon_{k})\cdot 2^k~(\hbox{\rm mod}~2^{k+1})$. But
${i_{k}+2^k\equiv 2^k-1+\varepsilon_{k}\cdot 2^k+2^k~(\hbox{\rm mod}~2^{k+1})}$. Thus 
$i_{k}+2^k\equiv -1+\varepsilon_{k}\cdot 2^k~(\hbox{\rm mod}~2^{k+1})$. Finally, 
$2^k-1\equiv i_{k}$ (resp., $i_{k}+2^k$) $(\hbox{\rm mod}~2^{k+1})$ 
if $\varepsilon_{k}=0$ (resp., $\varepsilon_{k}=1$).\bigskip 

\noindent - So $\varepsilon_{k}=0$ and $\gamma_{\alpha} (k)=\gamma_{\alpha'} 
(k)$.\bigskip

 This finishes the proof.\hfill{$\square$}\bigskip

 Now we prove that $[D_{2}(\boraone )\setminus \mbox{pot}(\bormone ),\leq^r_{c}]$ is not well-founded.\bigskip

\noindent {\bf Notation.}  ~Let 
${\cal S}:\omega^\omega\rightarrow\omega^\omega$ be the shift 
map: 
${\cal S}(\alpha )(k):=\alpha (k+1)$, $\beta_{0}$ be the 
sequence $(0,1,2,\ldots )$, and ${\beta_{n}:={\cal S}^n(\beta 
_{0})}$. Notice that $\beta_{n}(i)=i+n$, by induction on $n$. We put 
$B_{n}:=A^{S_{\beta_{n}}}$. 

\begin{prop} We have $B_{n+1}\leq^r_{c}B_{n}$ and 
$B_{n}\not\leq^r_{c}B_{n+1}$ for each integer $n$.\end{prop} 

\noindent {\bf Proof.} We define injective continuous maps 
$u=v:2^\omega\rightarrow 2^\omega$ by ${u(\alpha 
):=1^{1+n}\alpha}$. 
They are clearly witnesses for $B_{n+1}\leq^r_{c}B_{n}$.\bigskip

\noindent $\bullet$ Conversely, we argue by contradiction. This gives 
continuous maps $u$ and $v$ such that 
$${B_{n}=(u\times v)^{-1}(B_{n+1})}.$$ 
By Lemma 22, we have $u=v$. We set ${f^n_{m}:=f^{S_{\beta_{n}}}_{m}}$, and 
$f^n_{t}:=f^n_{t(0)}\ldots f^n_{t(|t|-1)}$ for 
$t\in\omega^{<\omega}\setminus\{\emptyset\}$, when it makes sense. 
Let 
$\alpha\in N_{0^{n+3}}$, so that ${\alpha = 0^{n+3}\gamma}$.\bigskip

\noindent $\bullet$ If $f^{n}_{t}(\alpha )$ is defined, then fix  
$m_{t}\!\in\!\omega$ with 
${u[f^{{n}}_{t}(\alpha )]\! =\! f^{{n+1}}_{m_{t}}
\big(u[f^{{n}}_{t-t\lceil 1}(\alpha )]\big)}$, and set ${U\!\! :=\! u(\alpha )}$ 
(with the convention that $f^{n}_{\emptyset}\! :=\!\mbox{Id}_{2^\omega}$). 
Then  
${u[f^{{n}}_{t}(\alpha )]\! =\! 
f^{{n+1}}_{m_{t},m_{t-t\lceil 1},\ldots ,m_{t-t\lceil (|t|-1)}}(U)}$. 
In particular, 
$$f^{{n+1}}_{m_{(1,\ldots ,n+2)},m_{(2,\ldots ,n+2)},\ldots  
,m_{n+2}}(U)=
f^{{n+1}}_{m_{(n+2,1,\ldots ,n+1)},m_{(1,\ldots ,n+1)},\ldots   
,m_{n+1}}(U).$$ 
Therefore
$$\{m_{(1,\ldots ,n+2)},m_{(2,\ldots ,n+2)},\ldots ,m_{n+2}\}=
\{m_{(n+2,1,\ldots ,n+1)},m_{(1,\ldots ,n+1)},\ldots ,m_{n+1}\}.$$ 

\vfill\eject

 If $m_{n+2}=m_{n+1}$, then we get $u(0^{n+2}1\gamma )=u(0^{n+1}10\gamma 
)$. As $f^{{n}}_{0}(0^{n+1}10\gamma )=10^n10\gamma$, we get 
$\big(u(0^{n+1}10\gamma ),v(10^n10\gamma )\big)\in B_{n+1}$ and 
$(0^{n+2}1\gamma ,10^n10\gamma )\in B_{n}$, which is absurd. Now 
suppose that 
$M:=\hbox{\rm max}(m_{(1,\ldots ,n+2)},m_{(2,\ldots ,n+2)},\ldots  
,m_{n+2})$ is in 
$\{m_{n+1},m_{n+2}\}$. Then we can find ${1\leq k\leq n+1}$ such that
$$\hbox{\rm Card}\big(U\lceil M\big),\ \hbox{\rm Card}\big(U\lceil M\big)\! +\! 
k\in\{\Sigma_{i<l}~\big( 1\! 
+\!\beta_{n+1}(i)\big)/l\!\in\!\omega\}.$$ 
But this is not possible, since 
${\Sigma_{i<l+1}~\big( 1\! +\!\beta_{n+1}(i)\big)\! 
-\!\Sigma_{i<l}~\big( 
1\! +\!\beta_{n+1}(i)\big)\! =\! l\! +\! n\! +\! 2}$.\bigskip

\noindent $\bullet$ We then get the contradiction by induction, since 
we can remove $M$ from both  
$$\{m_{(1,\ldots ,n+2)},m_{(2,\ldots ,n+2)},\ldots   ,m_{n+2}\}\mbox{,}$$ 
$\{m_{(n+2,1,\ldots ,n+1)},m_{(1,\ldots ,n+1)},\ldots ,m_{n+1}\}$.\hfill{$\square$}\bigskip

\noindent {\bf Remarks.}~(a) We showed that 
$(A^{S_{\beta_{\alpha}}})_{\alpha\in 2^\omega}$ is a perfect 
antichain made of 
sets minimal among non-$\mbox{pot}(\bormone )$ sets for $\leq^r_{c}$. There 
are other natural notions of reduction. 
We defined $\leq^r_{c}$ in the introduction. If we 
moreover ask that $u$ and $v$ are one-to-one, this defines a new 
quasi-order that 
we denote $\sqsubseteq^r_{c}$. If $u$ and $v$ are only Borel, we have 
two other 
quasi-orders, denoted $\leq^r_{B}$ and $\sqsubseteq^r_{B}$. If $X=Y$, 
$X'=Y'$ and $u=v$, we get the usual notions $\leq_{c}$, 
$\sqsubseteq_{c}$, $\leq_{B}$ and $\sqsubseteq_{B}$. Let $\leq$ be 
any of these eight quasi-orders. Then 
$(A^{S_{\beta_{\alpha}}})_{\alpha\in 2^\omega}$\bf ~is a perfect 
antichain made 
of sets minimal among non-$\mbox{pot}(\bormone )$ sets for $\leq$\rm 
:\bigskip 

\noindent $\bullet$ Let us go back to Theorem 15 first. Assume this 
time that 
$A\leq A^f$. Then in the second case we can have 
$A^f\cap B^2\sqsubseteq^r_{c}A$ if $\leq$ is rectangular, and 
$A^f\cap B^2\sqsubseteq_{c}A$ otherwise. The changes to make in the 
proof are the 
following. Let $\nu$ (resp., $\nu'$) be a finer Polish topology on 
$Y$ (resp., $Y'$) making $u$ (resp., $v$) continuous. We get continuous 
maps $u':2^\omega\rightarrow [Y,\nu ]$ and 
$v':2^\omega\rightarrow [Y',\nu']$. The proof shows that 
$f\vert_{G}$ and $g\vert_{G}$ are actually witnesses for $A^f\cap 
G^2\sqsubseteq^r_{c}A$ if $\leq$ is rectangular, and $A^f\cap G^2\sqsubseteq_{c}A$ 
otherwise.\bigskip

\noindent $\bullet$ In Corollary 16, we can replace $\leq^r_{c}$ with 
$\leq$.\bigskip

\noindent $\bullet$ The proof of Theorem 19 shows that, in its 
statement, we can write $A^{S}\sqsubseteq_{c} A^{S}\cap B^2$.\bigskip

\noindent $\bullet$ The proof of Lemma 20 shows that, in its statement,  
we can replace $\leq^r_{c}$ with $\leq$.\bigskip

\noindent $\bullet$ It follows from Corollary 21 that 
$A^{S_{\beta_{\alpha}}}$ 
is, in fact, minimal among non-$\mbox{pot}(\bormone )$ sets for 
$\leq$.\bigskip

\noindent $\bullet$ To see that 
$(A^{S_{\beta_{\alpha}}})_{\alpha\in 2^\omega}$ 
is an antichain for $\leq^r_{B}$, it is enough to see that in the 
statement of 
Theorem 23, we can replace $\perp^r_{c}$ with $\perp^r_{B}$. We only 
have to change the beginning of the proof of Theorem 23. This time $u$ and 
$v$ are 
Borel. Let $\tau$ be a finer Polish topology on $2^\omega$ making $u$ 
and $v$ continuous, and $X:=[2^\omega ,\tau ]$. By Lemma 20, $A^{S}$ is 
$\leq^r_{c}$-minimal, so $(2^\omega, A^{S})\leq^r_{c} (X, 
A^{S})\leq^r_{c} A^{S'}$,  
and we may assume that $u$ and $v$ are continuous.\bigskip

\noindent (b) We proved that 
$[D_{2}(\boraone )\setminus \mbox{pot}(\bormone ),\leq^r_{c}]$ is not 
well-founded. Let 
$\leq$ be any of the eight usual quasi-orders. Then 
$[D_{2}(\boraone )\setminus \mbox{pot}(\bormone ),\leq]$\bf ~is not 
well-founded\rm :\bigskip

\noindent $\bullet$ The proof of Proposition 25 shows that 
$B_{n+1}\sqsubseteq_{c}B_{n}$, thus $B_{n+1}\leq B_{n}$.\bigskip

\noindent $\bullet$ We have to see that $B_{n}\not\leq^r_{B} 
B_{n+1}$. We argue by 
contradiction, so that we get $u$ and $v$ Borel.\bigskip

\noindent $\bullet$ Let us show that we can find a dense $G_{\delta}$ 
subset $G$ of $2^\omega$ such that $u\vert_{G}=v\vert_{G}$ is continuous, and 
$f^{{n}}_{m}(\alpha )\in G$, for each $\alpha\in G\cap 
D_{f^{{n}}_{m}}$.\bigskip

\noindent {\bf Claim.}~\it The set ${H:=\{\alpha\!\in\! 2^\omega /\forall 
p~\exists m\!\geq\! p\ \ \alpha\!\in\! D_{f^{{n}}_{m}}\}}$ is a dense 
$G_{\delta}$ subset of $2^\omega$.\bigskip\rm 

 We argue by contradiction. We can find a nonempty clopen set $V$ 
disjoint from $H$. The set $B_{n}\cap V^2$ has finite sections, so is 
$\mbox{pot}(\bormone)$ (see Theorem 3.6 in [Lo1]). But $\big(V,({f^n_{m}}|_{V\cap 
{f^n_{m}}^{-1}(V)})\big)$ 
is a converging situation, so that $B_{n}\!\cap\! V^2$ is not 
$\mbox{pot}(\bormone)$.\hfill{$\diamond$}\bigskip 

 So we can find a dense $G_{\delta}$ subset $K$ of $2^\omega$ such that 
$u\vert_{K}$, $v\vert_{K}$ are continuous and $K\subseteq H$. 
Now let $K_{0}:=K$, 
$K_{p+1}:=K_{p}\setminus\big(\bigcup_{m}D_{f^{{n}}_{m}}\setminus 
{f^{{n}}_{m}}^{-1}(K_{p})\big)$, and 
${G:=\bigcap_{p}K_{p}}$. If $\alpha\in K_{1}$, fix $(m_{k})$  
infinite such 
that $\alpha\in\bigcap_{k}D_{f^{{n}}_{m_{k}}}$. We have 
$f^{{n}}_{m_{k}}(\alpha )\in K_{0}$, so $\big(u(\alpha 
),v[f^{{n}}_{m_{k}}
(\alpha )]\big)$ tends to ${\big(u(\alpha ),v(\alpha )\big)\in
\overline{B_{n+1}}\setminus B_{n+1}=\Delta (2^\omega )}$. So 
$u\vert_{K_{1}}=v\vert_{K_{1}}$. Now it is clear that $G$ is 
suitable.\bigskip 

\noindent $\bullet$ We take $\alpha\in G\cap N_{0^{n+3}}$ and  
complete the proof as we did for Proposition 25.\bigskip 

\noindent {\bf Proof~of~Theorem~6.} We will actually prove a stronger 
statement. We set
$$(P_{0},P_{1},P_{2},P_{3},P_{4}):=( \hbox{\rm 
reflexive,~irreflexive,~symmetric,~antisymmetric,~transitive}).$$
Let $\sigma\in 2^5\setminus\{\{2,4\},\{0,2,4\}\}$ such that the class 
$\Gamma$ of $\borel\setminus \mbox{pot}(\bormone)$ relations satisfying 
$\wedge_{j\in\sigma}~P_{j}$ 
is not empty. Then we can find a perfect $\leq_{B}$-antichain 
$(R_{\alpha})_{\alpha\in 2^\omega}$ in $D_{2}(\boraone)\cap\Gamma$ 
such that $R_{\alpha}$ is $\leq_{B}$-minimal among $\borel\setminus 
\mbox{pot}(\bormone)$ sets, for any $\alpha\in 2^\omega$.\bigskip

\noindent $\bullet$ First, notice that if $\{0,1\}\subseteq\sigma$ or 
$\sigma =\{1,2,4\}$, then every relation satisfying 
$\wedge_{j\in\sigma} P_{j}$ 
is empty, thus $\mbox{pot}(\bormone)$. If $\{2,3\}\subseteq\sigma$, then every 
Borel relation 
satisfying $\wedge_{j\in\sigma} P_{j}$ is a subset of the diagonal, and is 
therefore $\mbox{pot}(\bormone)$. If $\sigma =\{0,2,4\}$, we are in the case of Borel 
equivalence 
relations, and by Harrington, Kechris and Louveau's Theorem, $E_{0}$ 
is minimum among non-$\mbox{pot}(\bormone )$ equivalence relations. If 
$\sigma =\{2,4\}$, then any Borel relation $A\subseteq X^2$ 
satisfying $\wedge_{j\in\sigma} P_{j}$ is reflexive on its domain 
$\{x\in X/(x,x)\in A\}$, which is a Borel set. Thus we are reduced to 
the case of 
equivalence relations. In the sequel, we will avoid these cases and 
show the 
existence of a perfect antichain made of minimal sets for $[\Gamma 
,\leq_{B}]$.\bigskip

\noindent $\bullet$  Let $\mathbb{A}\!:=\!\{A^{S_{\beta_{\alpha}}}/\alpha\!\in\! 
2^\omega\}$. In the introduction, we defined $R_{A}$ for 
$A\subseteq 2^\omega\times 2^\omega$.\bigskip

\noindent {\bf Claim~1.}~\it $\{ R_{A}/A\!\in\! {\mathbb{A}}\}$ is a 
$\leq_B$-antichain.\bigskip\rm

 Assume that ${A\!\not=\! A'\!\in\! {\mathbb{A}}}$ satisfy 
$R_{A}\leq_B R_{A'}$. Then there is 
${f:2^\omega\times 2\rightarrow 2^\omega\!\times\! 2}$ with 
$${R_{A}=(f\!\times\! f)^{-1}(R_{A'}).}$$ 
We set $F_\varepsilon\! :=\!\{x\!\in\! 2^\omega\!\times\! 2/x_1\! =\!\varepsilon\}$ and 
${b_{\varepsilon}\! :=\! R_{A}\cap 
(F_\varepsilon\!\times\! F_\varepsilon )}$, for $\varepsilon\!\in\! 2$. 
We then put ${a\!:=\! R_{A}\cap (F_0\!\times\! F_1)}$. We have 
$R_{A}\! =\! a\cup b_{0}\cup b_{1}$, and 
${b_{\varepsilon}\! =\!\{(x,y)\!\in\! F_{\varepsilon}\!\times\! 
F_{\varepsilon}/x_{0}\! 
=\! y_{0}\}}$ is $\mbox{pot}(\bormone)$. We set 
$$F'_\varepsilon:=\{x\in 2^\omega\times 2/f_1(x)=\varepsilon\}$$ 
and ${b'_{\varepsilon}:=R_{A}\cap (F'_\varepsilon\times F'_\varepsilon )}$, for $\varepsilon\in 2$.

\vfill\eject

 We then put ${a':=R_{A}\cap (F'_0\times F'_1)}$. We have $R_{A}=a'\cup b'_{0}\cup b'_{1}$, and 
$${b'_{\varepsilon}=(f\vert_{F'_{\varepsilon}}\times f\vert_{F'_{\varepsilon}})^{-1}
\big(\Delta (2^\omega\times 2)\big)\in \mbox{pot}(\bormone)}.$$ 
Notice that $A\leq_{B}^r R_{A}$, so that $R_{A}$ is not $\mbox{pot}(\bormone)$. So 
${R_{A}=(a\cap a')\cup b_{0}\cup b_{1}\cup b'_{0}\cup b'_{1}}$, and 
$a\cap a'$ is not $\mbox{pot}(\bormone)$. It remains to define $C:=a\cap a'$, 
viewed as a subset of ${(F_0\cap F'_0)\times (F_1\cap F'_1)}$. We equip 
$F_0\cap F'_0$ 
(resp., $F_1\cap F'_1$) with a finer Polish topology making 
$f\vert_{F_0\cap F'_0}$ 
(resp., $f\vert_{F_1\cap F'_1}$) continuous. Then $C\!\leq^r_{c}\! 
A$ and $C\!\leq^r_{c}\! A'$, which contradicts Corollaries 21 and 
24.\hfill{$\diamond$}\bigskip

\noindent {\bf Claim~2.}~\it Let $A\! =\! A^{S_{\beta_{\alpha}}}\!\in\! 
{\mathbb{A}}$. 
Then $R_{A}$ is minimal for $\leq_{B}$ among $\borel\setminus 
\mbox{pot}(\bormone)$ relations.\bigskip\rm 

\noindent - Assume that $R\leq_{B} R_{A}$. This gives 
${f:X\rightarrow 2^\omega\!\times\! 2}$ Borel with 
${R\! =\! (f\!\times\! f)^{-1}(R_{A})}$. Again we set 
${F_{\varepsilon}:=\{x\in X/f_{1}(x)=\varepsilon\}}$ for 
$\varepsilon\in 2$, and 
we see that $R\cap (F_{0}\times F_{1})$ is not 
$\mbox{pot}(\bormone)$.\bigskip

\noindent - Let $\tau$ be a finer Polish topology on $X$ making $f$ 
continuous. By Theorem 9 there are 
$${u:2^\omega\rightarrow [F_{0},\tau ]\mbox{,}}$$ 
${v:2^\omega\rightarrow [F_{1},\tau ]}$ continuous with 
${A_{1}\! =\! (u\!\times\! v)^{-1}\big(R\cap (F_{0}\!\times\! 
F_{1})\big)\cap\overline{A_{1}}}$. 
We define ${H:=f_{0}\big[u[2^\omega ]\big]}$, 
${K:=f_{0}\big[v[2^\omega ]\big]}$ 
and $P:=H\cap K$; this defines compact subsets of $2^\omega$. Then 
${A\cap (H\times K)}$ is not $\mbox{pot}(\bormone)$ since
$${A_{1}\! =\! [(f_{0}\circ u)\!\times\! (f_{0}\circ 
v)]^{-1}\big(A\cap (H\times K)\big)\cap\overline{A_{1}}.}$$
As in the proof of Theorem 15, this implies that ${A\cap P^2}$ is not 
$\mbox{pot}(\bormone)$. By Lemma 14, we can find a Borel subset $S$ of $P$ 
and a finer topology $\sigma$ on $S$ such that $\big( [S,\sigma 
],(f^{S}_{n})_n\big)$ is a converging situation.\bigskip

\noindent - By 18.3 in [K], we can find a Baire measurable map 
${g_{\varepsilon}:S\rightarrow f^{-1}(S\times\{\varepsilon\})}$ 
such that 
$$f_{0}\big(g_{\varepsilon}(\alpha )\big)=\alpha\mbox{,}$$ 
for $\alpha$ in $S$ and 
$\varepsilon\in 2$. Let $G$ be a dense $G_{\delta}$ subset of $S$ 
such that each 
${g_{\varepsilon}\vert}_{G}$ is continuous. Now we define 
$F:G\times 2\rightarrow X$ by 
${F(\alpha ,\varepsilon ):=g_{\varepsilon}(\alpha )}$. Then  
$R_{A}\cap (G\times 2)^2=(F\times F)^{-1}(R)$, so 
$R_{A}\cap (G\times 2)^2\leq_{B}R$. As in the proof of Theorem 15, we 
see that $A\cap G^2$ is not $\mbox{pot}(\bormone)$. But $A\cap 
G^2\sqsubseteq_{c}A$. By Remark (a) above, we get $A\sqsubseteq_{c}A\cap G^2$. Thus 
$R_{A}\sqsubseteq_{c}R_{A}\cap (G\times 2)^2$ and ${R_{A}\leq_{B}R}$.
\hfill{$\diamond$}\bigskip

 Finally, one easily checks the existence of a continuous map 
$c:2^\omega\rightarrow 2^\omega$ such that $c(\delta )$ is a 
Borel code for 
$R_{A}$ if $\delta$ is a Borel code for $A$. So there is a continuous 
map ${r:2^\omega\rightarrow 2^\omega}$ such that $r(\alpha )$ is a 
Borel code for 
$R_{A^{S_{\beta_{\alpha}}}}$. This shows, in particular, the existence 
of a perfect antichain made of minimal sets for 
${[\borel\setminus \mbox{pot}(\bormone)~ \hbox{\rm quasi-orders},\leq_{B}]}$ and 
${[\borel\setminus \mbox{pot}(\bormone)~ \hbox{\rm partial~orders},\leq_{B}]}$. More 
generally, this works if $\sigma\subseteq\{0,3,4\}$.\bigskip

\noindent $\bullet$ Similarly, we define, for $A\subseteq X^2$, a 
strict partial order relation $R'_A$ on $X\times 2$ by
$$(x,i)~R'_A~(y,j)~~\Leftrightarrow ~~[(x,y)\in A~ \hbox{\rm and}~
i=0~ \hbox{\rm and}~j=1].$$
The proof of the previous point shows that if 
$\sigma\subseteq\{1,3,4\}$, 
then $\{R'_{A}/A\in {\mathbb{A}}\}$ is a perfect antichain made of 
minimal sets for 
$[\Gamma ,\leq_{B}]$. Notice that this applies when $\Gamma$ is 
the class of $\borel\setminus \mbox{pot}(\bormone)$ strict quasi-orders, 
strict partial orders, directed graphs or oriented graphs.\bigskip

\noindent $\bullet$ Similarly again, we can define, for $A\subseteq 
X^2$, $S_A$  
reflexive symmetric on $X\times 2$ by
$$\begin{array}{ll}  
(x,i)~S_A~(y,j)~~\Leftrightarrow 
& (x,i)=(y,j)~ \hbox{\rm or}~
[(x,y)\in A~ \hbox{\rm and}~i=0~ \hbox{\rm and}~j=1]~ \hbox{\rm 
or}~\cr 
& [(y,x)\in A~ \hbox{\rm and}~i=1~ \hbox{\rm and}~j=0]. 
\end{array}$$
Let $A_{0}:=A$ and $A_{1}:={A}^{-1}$. The proof of Claim 1 shows that 
if $A\not= A'\in {\mathbb{A}}$ satisfy $S_{A}\leq_B S_{A'}$, then we 
can find $C\notin \mbox{pot}(\bormone )$ and $\varepsilon ,~\varepsilon'\in 
2$ such that $C\leq^r_{c} A_{\varepsilon}$ and $C\leq^r_{c} 
A'_{\varepsilon'}$. But this contradicts Corollaries 21 and 24. This shows that 
if $\sigma=\{0,2\}$, then 
$\{S_{A}/A\in {\mathbb{A}}\}$ is a perfect antichain made of minimal sets 
for $[\Gamma ,\leq_{B}]$.\bigskip

\noindent $\bullet$ Similarly again, we can define, for $A\subseteq 
X^2$, 
a graph relation $S'_A$ on $X\times 2$ by  
$$(x,i)~S'_A~(y,j)~\Leftrightarrow ~
[(x,y)\!\in\! A~ \hbox{\rm and}~i\! =\! 0~ \hbox{\rm and}~j\! =\! 1]~ 
\hbox{\rm or}~ 
[(y,x)\!\in\! A~ \hbox{\rm and}~i\! =\! 1~ \hbox{\rm and}~j\! =\! 0].$$
The proof of the previous point shows that if 
$\sigma\subseteq\{1,2\}$, then 
$\{S'_{A}/A\in {\mathbb{A}}\}$ is a perfect antichain made of minimal 
sets for 
$[\Gamma ,\leq_{B}]$. Notice that this applies when $\Gamma$ is the class of 
$\borel\setminus \mbox{pot}(\bormone)$ graphs. This finishes the proof.\hfill{$\square$}\bigskip

\noindent {\bf Remarks.}~(a) We showed that $(R_{A})_{A\in {\mathbb{A}}}$ 
is a perfect 
antichain made of sets minimal among non-$\mbox{pot}(\bormone )$ sets for 
$\leq_{B}$. Fix $\leq$ in 
$\{\leq_{c},\sqsubseteq_{c},\leq_{B},\sqsubseteq_{B}\}$. Then 
$(R_{A})_{A\in {\mathbb{A}}}$\bf\ is a perfect antichain made of sets 
minimal among 
non-$\mbox{pot}(\bormone )$ sets for $\leq$.\rm\ It is enough to check the 
minimality. The 
only thing to notice, in the proof of Claim 2 of the proof of Theorem 
6, is that we have $R_{A}\cap (G\times 2)^{2}\sqsubseteq_{c}R$ and 
$R_{A}\sqsubseteq_{c} R$. Similarly, 
$R'_{A}$, $S^{}_{A}$ and $S'_{A}$ ($A\in {\mathbb{A}}$) are minimal for 
$\leq_{c}$, $\sqsubseteq_{c}$ and $\sqsubseteq_{B}$.\bigskip

\noindent (b) We have $\neg\Delta (2^\omega )\perp^r_{B}L_{0}$. 
Indeed, assume that $\neg\Delta (2^\omega )=(u\times v)^{-1}(L_0)$. Then 
$u(\alpha )<_{\hbox{\rm lex}} v(\beta )$ if 
$\alpha\not=\beta$, and $v(\alpha )\leq_{\hbox{\rm lex}} u(\alpha )$. Thus 
$$u(\beta )<_{\hbox{\rm lex}} v(\alpha )\leq_{\hbox{\rm lex}} u(\alpha )
<_{\hbox{\rm lex}}v(\beta )\leq_{\hbox{\rm lex}} u(\beta )\mbox{,}$$ 
which is absurd. Now assume that 
$L_0=(u\times v)^{-1}\big(\neg\Delta (2^\omega )\big)$. Then 
$\beta\leq_{\hbox{\rm lex}}\alpha$ implies $u(\alpha ) = v(\beta )$, thus $u=v$ 
has to be 
constant. Thus $\alpha<_{\hbox{\rm lex}}\beta$ implies that $u(\alpha )$ and 
$v(\beta )$ are different and equal.\bigskip

 In the introduction, we saw that $\{\neg\Delta (2^\omega ),L_{0}\}$ 
is a complete family of minimal sets for 
$${[\mbox{pot}\big(\check 
D_2(\boraone )\big)\setminus \mbox{pot}(\bormone ),\sqsubseteq^r_{c}].}$$ 
We just saw that $\{\neg\Delta (2^\omega ),L_{0}\}$ is an antichain for 
$\leq^r_{B}$, and 
therefore for any of the eight usual quasi-orders. These facts imply 
that $\neg\Delta (2^\omega )$ and $L_{0}$ are minimal among 
non-$\mbox{pot}(\bormone )$ sets 
for $\leq^r_{c}$, $\sqsubseteq^r_{c}$, $\leq^r_{B}$ and 
$\sqsubseteq^r_{B}$. But 
$\neg\Delta (2^\omega )$ and $L_{0}$ are also minimal for $\leq_{c}$, 
$\sqsubseteq_{c}$, $\leq_{B}$ and $\sqsubseteq_{B}$. Indeed, if $O$ 
is any of 
these two open sets, we have $\overline{O}\setminus O=\Delta 
(2^\omega )$. This 
gives $G$ such that $O\cap G^2\sqsubseteq_{c}A$, as in the proof of 
Theorem 15 (and Remark (a) after Proposition 25). Then any increasing 
continuous injection $\phi: 2^\omega\rightarrow G$ is a 
witness to $O\sqsubseteq_{c}O\cap G^{2}$.\bigskip\smallskip

\noindent\bf {\Large 5 The minimality of $A_{1}$ for the classical notions of 
comparison.}\rm\bigskip

 As announced in the introduction, we will show a result implying 
that $A_{1}$ is minimal among non-potentially closed sets. The following 
definition specifies the meaning of the expression 
``the $f_{n}$'s do not induce cycles" mentioned in the 
introduction. This kind of notion has already been used in the theory 
of potential complexity (see Definition 2.10 in [L3]).  

\begin{defi}  We say that $\big(X,(f_{n})\big)$ is an $acyclic\ situation$ if\smallskip

\noindent (a) $\big(X,(f_{n})\big)$ is a converging situation, with only 
$\Delta (X)\subseteq\overline{A^f}\setminus A^f$ in condition (c).\smallskip  

\noindent (b) For $v\in\omega^{<\omega}\setminus\{\emptyset\}$ and 
$\varepsilon\in\{-1,1\}^{|v|}$, the following implication 
holds:\smallskip
  
\leftline{$\big(\forall i\! <\! |v|\! -\! 1~~v(i)\!\not=\! v(i\! +\! 1)\ 
\mbox{or}\ \varepsilon (i)\!\not=\! -\varepsilon (i\! +\! 1)\big)
\Rightarrow 
\big(\forall U\!\in\!\borone (X)\!\setminus\!\{\emptyset\}~\exists 
V\!\in\!\borone (U)\!\setminus\!\{\emptyset\}$}\smallskip
 
\rightline{$\forall x\!\in\! V\ [f_{v(|v|-1)}^{\varepsilon (|v|-1)}\ldots f_{v(0)}^{\varepsilon (0)}
(x)~\mbox{is~not~defined~or~not~in}~V ]\big).$}\end{defi}

\noindent\bf Notation.\rm\ We define 
$f^1_{n}:N_{s_{n}0}\rightarrow N_{s_{n}1}$ by 
$f^1_{n}(s_{n}0\gamma ):=s_{n}1\gamma$ (where $s_n$ is as defined in the 
introduction, to build $A_1=\bigcup_n\ \hbox{Gr}(f^1_n)$).

\begin{lem} Let $\alpha\in 2^\omega$, 
$v\in\omega^{<\omega}\setminus\{\emptyset\}$ and 
$\varepsilon\in\{-1,1\}^{|v|}$. Assume that ${v(i)\not= v(i+1)}$ or 
$\varepsilon (i)\not= -\varepsilon (i+1)$ if $i<|v|-1$. Then 
${f^{1}_{v(|v|-1)}}^{\varepsilon (|v|-1)}\ldots {f^{1}_{v(0)}}^{\varepsilon (0)}
(\alpha )$ is either undefined, or of value different than $\alpha$.
\end{lem}

\noindent\bf Proof.\rm\ We argue by contradiction. Let $v$ be a 
counter-example of minimal length. Note that $|v|\geq 3$. Set 
$l:=\hbox{\rm max}_{i<|v|}~v(i)$, $e_0:=e_{|v|}:=\alpha\lceil (l\! +\! 1)$, and, 
for $0<i<|v|$: 
$$e_i:=[{f^{1}_{v(i-1)}}^{\varepsilon (i-1)}\ldots {f^{1}_{v(0)}}^
{\varepsilon (0)}(\alpha )]~\lceil ~(l\! +\! 1).$$ 
Set $\Theta :=(\theta_n)$, where $\theta_n:=s_n$. Then $(e_i)_{i\leq |v|}$ is an 
$s({\mathfrak R}_\Theta )$-cycle, which contradicts Proposition 
18.(b).\hfill{$\square$}\bigskip

\noindent\bf Example.\rm\ $\big(2^\omega ,(f^1_{n})\big)$ 
is an acyclic situation. Indeed, $\big(2^\omega ,(f^{1}_n)\big)$ is a converging 
situation, by Corollary 12. Let 
us show that condition (b) in the definition of an acyclic situation 
is true for $\big(2^\omega ,(f^{1}_n)\big)$. The domain $D$ of 
${f^{1}_{v(|v|-1)}}^{\varepsilon (|v|-1)}\!\!\!\!\ldots 
{f^{1}_{v(0)}}^{\varepsilon (0)}$ 
is clopen. If $U$ is not included in $D$, then we can take $V:=U\setminus 
D$. Otherwise, 
let $\alpha\!\in\! U$. By Lemma 27, and by continuity, we can find a 
clopen neighborhood $V$ of $\alpha$ included in $U$ such that 
${f^{1}_{v(|v|-1)}}^{\varepsilon (|v|-1)}\ldots 
{f^{1}_{v(0)}}^{\varepsilon (0)}[V]\cap V=\emptyset$.

\begin{thm} Let $\big(X,(f_{n})\big)$ be an acyclic 
situation. Then $A_{1}\leq^r_{c}A^f$.\end{thm}

\noindent {\bf Proof.} It looks like those of Theorems 2.6 and 2.12 
in [L3]. The main difference is that we want a reduction defined on the whole 
product. It is also similar to the proof of Theorem 19. Let us indicate the 
differences with the proof of Theorem 19. We replace $A^S=\bigcup_n\ 
\hbox{Gr}(f^S_n)$ with $A_1=\bigcup_n\ \hbox{Gr}(f^1_n)$.\bigskip

\noindent $\bullet$ We only construct $(U_{s})_{s\in 2^{<\omega}}$ 
and $\phi$, so that (iii) becomes
$$(iii)\ (s~{\mathfrak  R}~t~\hbox{\rm and}~s\!\not=\! t)\Rightarrow 
U_{t}\! =\! f_{\phi (|s\wedge t|)}[U_{s}].$$
\noindent $\bullet$ Here we choose $\Theta =(\theta_n)$ with 
$\theta_n:=s_n$. Notice that 
${\mathfrak R}_\Theta ={\mathfrak R}$.\bigskip

\noindent $\bullet$ Condition $(3)$ becomes
$$(3)\ (t_{k}~{\mathfrak R}~t_{l}~\hbox{\rm and}~t_{k}\!\not=\! 
t_{l})\Rightarrow 
U^n_{l}\! =\! f_{\phi (|t_{k}\wedge t_{l}|)}[U^n_{k}].$$

\vfill\eject

\noindent $\bullet$ We can find $C\in\borone (U_{t_{0}\lceil 
p})\setminus\{\emptyset\}$ such that 
$C^{2}\cap\bigcup_{q\leq p}~\hbox{Gr}(f_{q})=\emptyset$, and also 
$${{n_{0}}\geq\hbox{\rm sup}~\{\phi (q)\! +\! 1/q<p\}}$$ 
with ${C^{2}\cap \hbox{Gr}(f_{n_{0}})\!\not=\!\emptyset}$, since  
${\Delta (X)\subseteq\overline{A^f}\setminus A^f}$. We set 
$\phi (p):={n_{0}}$. We first construct clopen sets $\tilde U^n_k$ as in 
the proof of Theorem 19.\bigskip

\noindent {\bf Case 2.}~ $t_{n}\lceil p\not=t_{r}\lceil p$.\bigskip

\noindent 2.1. $t_{r}~{\mathfrak R}_\Theta ~t_{n}$.\bigskip

 To get condition (4), fix $k,l\!\leq\! n$ with $\neg ~t_{k}~{\mathfrak R}~t_{l}$. 
Set ${f^\varepsilon_{v}\!:=\! {f_{\phi (v(|v|-1))}^{\varepsilon (|v|-1)}\ldots 
f_{\phi (v(0))}^{\varepsilon (0)}}}$, so that 
${\tilde U^n_{l}={f^{\varepsilon}_{v}}[\tilde U^n_{k}]}$, and we have 
$\phi\big( v(i)\big)\not=\phi\big( v(i+1)\big)$, since $\phi$ is 
strictly increasing. As 
$\big(X,(f_{n})\big)$ is without cycles, we can find $x\in\tilde 
U^n_k$ with 
$f^{\varepsilon}_{v}(x)\not= x$. We can therefore find a clopen 
neighborhood $\underline{U}^n_k$ of $x$, included in $\tilde U^n_k$, 
such that $\underline{U}^n_k\cap f^{\varepsilon}_{v}[\underline{U}^n_k] 
=\emptyset$. We construct clopen sets $\underline{U}^n_r$, for $k\not= r\leq n$, 
as before, ensuring condition (3). Notice that $\underline{U}^n_r\subseteq\tilde 
U^n_r$, so that the hereditary conditions (1), (2) and (5) remain 
fullfilled. In finitely many steps we get 
$(\underline{U}^n_k\times\underline{U}^n_l)\cap\Delta (X)=\emptyset$, 
for each pair $(k,l)$. The argument is similar for $\hbox{Gr}(f_{q})$ instead of 
$\Delta (X)$.\bigskip 

\noindent 2.2. $t_{n}~{\mathfrak R}_\Theta ~t_{r}$.\bigskip

 This case is similar to case 2.1.\hfill{$\square$}\bigskip

\noindent {\bf Remark.}~We actually showed that 
$A_{1}\sqsubseteq_{c}A^f$.

\begin{cor} $A_{1}$ is minimal among non-potentially 
closed sets for the eight usual quasi-orders.\end{cor}

\noindent {\bf Proof.} Let $B\!\in\!\borel (2^\omega )$, $\tau$ a 
finer topology on $B$, $Z\!:=\! [B,\tau ]$ and 
${f_{n}\!:=\! {f^{1}_{n}}\vert_{B\cap {f^{1}_{n}}^{-1}(B)}}$. We 
assume that $\big(Z,(f_{n})\big)$ is a 
converging situation. By Corollary 16 and Remark (a) after Proposition 
25, it is enough to show that ${A_{1}\sqsubseteq_{c}A_{1}\cap Z^{2}=A^f}$. By 
Theorem 28 and the remark above, it is enough to check that $\big(Z,(f_{n})\big)$ 
is an acyclic situation, i.e., condition (b). Fix $\alpha\!\in\! U$ and 
$f_{v}^\varepsilon:=f_{v(|v|-1)}^{\varepsilon (|v|-1)}\ldots 
f_{v(0)}^{\varepsilon (0)}$. If $U$ is not included in 
$D_{f_{v}^\varepsilon}$, then we 
can take $V:=U\setminus D_{f_{v}^\varepsilon}$, because the domain is 
a clopen subset of $Z$. As $f_{v}^\varepsilon$ is continuous, it is enough to 
see that $f_{v}^\varepsilon (\alpha )\not=\alpha$, if $U$ is included in 
$D_{f_{v}^\varepsilon}$. But this is clear, since 
${f^{1}_{v(|v|-1)}}^{\varepsilon (|v|-1)}\ldots
{f^{1}_{v(0)}}^{\varepsilon (0)}(\alpha )$ is different from 
$\alpha$, by Lemma 27.\hfill{$\square$}\bigskip

\noindent {\bf Remarks.}~(a) Theorem 28 is also a consequence of the 
following result:

\begin{thm} (Miller) Let $X$ be a Polish space, and $A$ a locally 
countable $\ana$ oriented graph on $X$ whose symmetrization is acyclic (in the 
sense of Definition 17). Then exactly one of the following holds:\smallskip

\noindent (a) $A$ has countable Borel chromatic number.\smallskip

\noindent (b) $A_{1}\sqsubseteq_{c}A$.\end{thm}

 Theorem 30 is actually a corollary of a more general result, motivated by the 
results of this paper, which gives a basis for locally countable Borel directed 
graphs of uncountable Borel chromatic number, with respect to $\sqsubseteq_{c}$. 
The proof of both Theorem 30 and the basis result appear in [M1].

\vfill\eject\bigskip

\noindent (b) We saw that $A_{1}\sqsubseteq_{c}A^f$ if $\big(X,(f_{n})\big)$ is 
an acyclic situation. There is another example of a 
$$D_{2}(\boraone )\setminus \mbox{pot}(\bormone )$$ 
set, which seems more ``natural" than $A_{1}$. It is
$$C_{1}\! :=\!\{(\alpha ,\beta )\!\in\! 2^\omega\!\times\! 2^\omega 
/\exists s\!\in\! 2^{<\omega}~\exists\gamma\!\in\! 2^\omega ~~
(\alpha ,\beta )\! =\! (s0\gamma ,s1\gamma )\}.$$
Its symmetric version plays an important role in the theory of potential 
complexity (see for example Theorem 3.7 and Corollary 4.14 in [L1]). We wonder 
what $\{ C_{1}\}$ is a basis for. Roughly speaking, $\{ C_{1}\}$ will be a basis 
for situations where commuting relations between the $f_{n}$'s are involved. More 
specifically, 

\begin{defi} We say that $\big(X,(f_{n})\big)$ 
is a $commuting~situation$ if\smallskip

\noindent (a) $X$ is a nonempty perfect closed subset of $\omega^\omega$.\smallskip

\noindent (b) $f_{n}$ is a partial homeomorphism with disjoint $\borone(X)$ 
domain and range. Moreover $\alpha <_{\mbox{lex}} f_{n}(\alpha )$ if 
$\alpha\in D_{f_{n}}$.\smallskip

\noindent (c) $\Delta (X)\subseteq\overline{A^f}\setminus A^f$, and 
$A^f\in\bormtwo (X^{2})$.\smallskip
  
\noindent (d) For each $\alpha\in f_{m}^{-1}(D_{f_{n}})$ we have 
$\alpha\in f_{n}^{-1}(D_{f_{m}})$ and 
$f_{m}\big(f_{n}(\alpha )\big)=f_{n}\big(f_{m}(\alpha )\big)$. Moreover the 
graphs of the $f_{n}$'s are pairwise disjoint.\end{defi}

 A $0$-dimensional Polish space is homeomorphic to a closed subset of 
$\omega^\omega$. So condition (a) is essentially the same as condition 
(a) of a converging situation. We use this formulation for the last part of  
condition (b). The disjunction of the domain and the range of $f_{n}$, and the 
inequality $\alpha <_{\hbox{\rm lex}} f_{n}(\alpha )$ come from symmetry 
problems. We will 
come back later to this. We will also come back to the $\bormtwo$ condition. It 
is linked with transitivity properties. The first part of condition (d) 
expresses the commutativity of the functions. One can show the following 
result, whose proof contains a part quite similar to the proof of Theorems 19 and 
28. 

\begin{thm} Let $\big(X,(f_{n})\big)$ be a commuting 
situation. Then $C_{1}\sqsubseteq_{c}A^f$.\end{thm}

 The proof of this uses the fact that $C_{1}=A^{f}$, where 
$\big(2^\omega ,(f_{n})\big)$ is a commuting situation. Let 
$$g_{n}:2^\omega\rightarrow 2^\omega$$ 
be defined by $g_{n}(\alpha )(k):=\alpha (k)$ if $k\not= n$, $1-\alpha (n)$ otherwise. Then 
$s(C_{1})=\bigcup_{n}~\hbox{\rm Gr}(g_{n})$, so $\big(2^\omega ,(g_{n})\big)$ is 
not a commuting situation, since otherwise we would have 
$C_{1}\sqsubseteq_{C} s(C_{1})$, which is absurd since $s(C_{1})$ is symmetric 
and $C_{1}$ is not. But the two reasons for that are that 
$\alpha\not <_{\hbox{\rm lex}} g_{n}(\alpha )$, and that the domain and the range 
of the bijections $g_{n}$ are not disjoint.\bigskip

 Similarly, let $\phi :\omega\rightarrow P_{f}\setminus\{ 0^\infty\}$ be a 
bijective map. We let $g'_{n}(\alpha )(p)\! :=\!\alpha (p)$ if 
$\phi (n)(p)\! =\! 0$, $1$ otherwise. This defines 
${g'_{n}\! :\!\{\alpha\!\in\! 2^\omega\! /\forall p~\phi (n)(p)\! =\! 0~or~
\alpha (p)\! =\! 0\}\!\rightarrow\! 2^\omega\! }$. Note that  
$${E_{0}\!\cap\! L'_{0}\! =\!\bigcup_{q} \hbox{\rm Gr}(g'_{n})\mbox{,}}$$ 
where $L'_{0}\! :=\!\{ (\alpha ,\beta )\!\in\! 2^\omega\!\times\! 2^\omega\! /
\forall i\!\in\!\omega ~\alpha (i)\!\leq\!\beta 
(i)~\mbox{and}~\alpha\!\not=\!\beta\}$.

\vfill\eject
 
 Then $\big(2^\omega ,(g'_{n})\big)$ is not a commuting situation, since otherwise 
$C_{1}\sqsubseteq_{C} E_{0}\cap L'_{0}$, which is absurd since 
$E_{0}\cap L'_{0}$ is transitive and $C_{1}$ is not. But the reason for that is 
that $E_{0}\cap L'_{0}\notin\bormtwo$.\bigskip 

 B. D. Miller has also a version of Theorem 32 for directed graphs of uncountable 
Borel chromatic number (in [M2]). Its proof uses some methods analogous to those 
in the proof of Theorem 30. All of this shows the existence of 
numerous analogies between non potentially closed directed graphs and 
directed graphs of uncountable Borel chromatic number.\bigskip\bigskip 

\noindent\bf {\Large 6 References.}\rm\bigskip

\noindent [H-K-Lo]\ \ L. A. Harrington, A. S. Kechris and A. Louveau,~\it A Glimm-Effros 
dichotomy for Borel equivalence relations,~\rm J. Amer. Math. Soc.~3 (1990), 903-928

\noindent [K]\ \ A. S. Kechris,~\it Classical Descriptive Set Theory,~\rm 
Springer-Verlag, 1995

\noindent [K-S-T]\ \ A. S. Kechris, S. Solecki and S. Todor\v cevi\'c,~\it Borel chromatic numbers,\ \rm 
Adv. Math.~141 (1999), 1-44

\noindent [L1]\ \ D. Lecomte,~\it Classes de Wadge potentielles et th\'eor\`emes d'uniformisation 
partielle,~\rm Fund. Math.~143 (1993), 231-258

\noindent [L2]\ \ D. Lecomte,~\it Uniformisations partielles et crit\`eres \`a la Hurewicz dans le plan,~\rm Trans. A.M.S. ~347, 11 (1995), 4433-4460

\noindent [L3]\ \ D. Lecomte,~\it Tests \`a la Hurewicz dans le plan,~\rm Fund. Math.~156 (1998), 
131-165

\noindent [L4]\ \ D. Lecomte,~\it Complexit\'e des bor\'eliens~\`a coupes 
d\'enombrables,~\rm Fund. Math.~165 (2000), 139-174

\noindent [Lo1]\ \ A. Louveau,~\it A separation theorem for $\Ana$ sets,~\rm Trans. A. M. S.~260 (1980), 363-378

\noindent [Lo2]\ \ A. Louveau,~\it Ensembles analytiques et bor\'eliens dans les espaces produit,~\rm 
Ast\'erisque (S. M. F.) 78 (1980)

\noindent [Lo-SR]\ \ A. Louveau and J. Saint Raymond,~\it Borel classes and closed games : Wadge-type and Hurewicz-type results,~\rm Trans. A. M. S.~304 (1987), 431-467

\noindent [M1]\ \ B. D. Miller,~\it Basis theorems for graphs of uncountable Borel chromatic number,~\rm pre-print

\noindent [M2]\ \ B. D. Miller,~\it A dichotomy theorem for graphs induced by commuting families 
of Borel injections,~\rm pre-print

\noindent [M]\ \ Y. N. Moschovakis,~\it Descriptive set theory,~\rm North-Holland, 1980

\noindent [SR]\ \ J. Saint Raymond,~\it La structure bor\'elienne d'Effros est-elle standard ?,~\rm 
Fund. Math.~100 (1978), 201-210

\bigskip\bigskip

\noindent\bf Acknowledgements.\rm ~I would like to thank A. Louveau 
for his interest for this work, and also for some nice remarks that 
can be found in this paper. I thank B. D. Miller who helped me to 
improve the English in this paper while I was visiting UCLA during June 
2005. I also thank B. D. Miller for some nice remarks that improved the quality 
of this article.

\end{document}